\documentclass[a4paper,11pt]{article}
\usepackage[a4paper,left=3cm,right=3cm,top=3cm,bottom=3cm]{geometry}
\usepackage{amsmath, amssymb, amsthm}
\numberwithin{equation}{section} 

\usepackage{newtxtext,newtxmath}
\usepackage{mathtools}
\usepackage{physics}
\usepackage{esint}
\usepackage{subcaption}
\usepackage{graphicx}
\usepackage{enumitem}
\usepackage{xcolor}
\usepackage{hyperref}
\usepackage{cleveref}
\usepackage{comment}
\usepackage{orcidlink}
\usepackage{needspace}

\title{Spectral bounds for the operator pencil\\ of an elliptic system in an angle}
\author{
Michael Tsopanopoulos\, \orcidlink{0009-0007-3167-6862}\thanks{Weierstrass Institute for Applied Analysis and Stochastics, Anton-Wilhelm-Amo-Str. 39, 10117 Berlin, Germany.\\
\texttt{tsopanopoulos@wias-berlin.de}}
}
\date{April 28, 2025}

\newtheorem{theorem}{Theorem}[section]
\newtheorem{lemma}[theorem]{Lemma}
\newtheorem{proposition}[theorem]{Proposition}

\theoremstyle{definition}
\newtheorem{definition}[theorem]{Definition}
\newtheorem{example}[theorem]{Example}
\theoremstyle{remark}
\newtheorem*{remark}{\textbf{Remark}}

\newcommand{\Mat}{\operatorname{Mat}}
\newcommand{\Id}{\operatorname{Id}}

\begin{document}

\maketitle

\begin{abstract}
The model problem of a plane angle for a second-order elliptic system subject to Dirichlet, mixed, and Neumann boundary conditions is analyzed. For each boundary condition, the existence of solutions of the form $r^\lambda v$ is reduced to spectral analysis of a particular matrix. Focusing on Dirichlet and mixed boundary conditions, optimal bounds on $|\Re \lambda|$ are derived, employing tools from numerical range analysis and accretive operator theory. The developed framework is novel and recovers known bounds for Dirichlet boundary conditions. The results for mixed boundary conditions are new and represent the central contribution of this work. Immediate applications of these findings are new regularity results for linear second-order elliptic systems subject to mixed boundary conditions.
\end{abstract}

numerical range, accretive operators

\medskip


\section{Introduction}
Regularity theory for partial differential equations is concerned with the question of how regular a solution can be with respect to the input data, such as the source function and the boundary data. A classic example is the Laplace equation, $\Delta u=f$, on a domain $\Omega$ with smooth boundary $\partial \Omega$, subject to either Dirichlet or Neumann boundary conditions. In this setting, it is established that $u\in W^{k+2,p}(\Omega)$ when $f\in W^{k,p}(\Omega)$, for any $k\in \mathbb{N}$ and $p>1$ (§2.4 in \cite{GrisEPND2011}). However, as demonstrated in \cite{Sham1968RMEP}, this result fails in the presence of mixed boundary conditions. Also, for pure Dirichlet or Neumann boundary conditions, the results fails for polygonal domains $\Omega$ that include edges or vertices. Interestingly, the breakdown in regularity across these scenarios exhibits a common structure: Consider a cone $\mathcal{K} \subset \mathbb{R}^n$ subject to Dirichlet, Neumann, or mixed boundary conditions. Near the vertex of $\mathcal{K}$, a solution $u$ to a linear elliptic system  can be asymptotically described by terms of the form (neglecting factors of $\log r$, see \cite{Daug1988EBVP}, \cite{GrisEPND2011}, \cite{KMR1997EBVP}):
\begin{align}\label{asm}
    u(r,\omega) &\sim r^\lambda v(\omega), \\ \nonumber
    \text{where }(r, \omega) &\text{ are spherical coordinates with } r \text{ being the distance to the vertex.}
\end{align}
Here, $v$ is a function defined on a subset of the sphere $\mathbb{S}^{n-1}$ (the cone opening), and $\lambda\in \mathbb{C}$ is the crucial regularity parameter - as $|\Re \lambda|$ increases, the integrability and differentiability of $u$ improves. To this model problem one can associate the so-called operator pencil $\mathcal{A}(\lambda)$, a $\lambda$-dependent differential operator satisfying $\mathcal{A}(\lambda)v=0$, i.e., the spectrum of $\mathcal{A}$ consists of all possible exponents $\lambda$ in the asymptotic expansion (\ref{asm}). Consequently, bounds on the eigenvalues of $\mathcal{A}(\lambda)$ yield information about the regularity of $u$. By localization, the regularity of boundary value problems in general polyhedral domains can be reduced to model problems in cones and angles. In other words, the (lack of) regularity of solutions to elliptic systems is characterized by the asymptotic behavior in (\ref{asm}), or equivalently, by the spectrum of the corresponding operator pencil $\mathcal{A}$. A well-established reference in this regard is \cite{VMR2001SPCS}, which provides estimates for $\Re \lambda$ for various model problems. The follow-up work \cite{MaRO2010EEPD} applies these results to specific problems in three-dimensional polyhedral domains, translating the estimates for $\Re \lambda$ into regularity results.\medskip

The primary motivation for this paper was to improve the known regularity estimates for solutions in linear elasticity. From the considerations above, it is evident that this is inherently linked to solutions of the form $r^\lambda v$, arising from model problems that reflect the geometry of the original problem. The underlying equations can be delicate because elasticity involves systems of elliptic equations, as the displacement is a vector field, and typically incorporates mixed boundary conditions. Existing literature on regularity theory for linear elasticity is sparse, with many results providing only relatively weak estimates or lacking full generality. The regularity of solutions $u$ for the special case of the Lamé system is discussed in \cite{Nica1992LSPP} for three-dimensional polygonal domains. Here, geometrical conditions are found to ensure $u\in W^{1,p}$ for some $p>3$. In \cite{HMWW2019HRSE}, it is shown for arbitrary space dimensions and general elliptic systems of second order (without using the expansion $r^\lambda v$ at all) that $u\in W^{1,p}$ for some $p>2$, which is only a slight improvement with respect to $p=2$. In \cite{HDKREMPM2008}, results on three-dimensional \textit{scalar} elliptic model problems are given, including mixed boundary conditions, that yield $p>3$. In contrast, the counter-example in \cite{Sham1968RMEP} demonstrates an upper bound, showing that one cannot expect $p\geq 4$ for scalar equations in the two-dimensional half-space, and consequently not for elliptic systems.\medskip

This work develops a framework for studying the model problem in a two-dimensional angle for second-order elliptic systems with real-valued coefficient matrices. Utilizing this framework, we recover well-known bounds for Dirichlet boundary conditions in a novel way. The main result of this work, however, focuses on the model problem with mixed boundary conditions (Theorem \ref{main mixed}). For this case, we prove - under mild ellipticity conditions - that any solution of the form $r^\lambda v$ to the model problem satisfies the bounds $ |\Re \lambda|\geq \frac{1}{2}$ for $\alpha\leq \pi$ and $| \Re \lambda|\geq \frac{1}{4}$ for $\alpha\leq 2\pi$, where $\alpha$ is the opening angle. These bounds translate to regularity results in two-dimensional domains of the form $u\in W^{1,p}$ for $p=4-\varepsilon$, resp. $p=\frac{8}{3}-\varepsilon$. Moreover, these results can be applied to three-dimensional domains in which a change of boundary condition occurs along an edge, since the geometry reduces to model problems in angles. For brevity, we outline a simple scenario from linear elasticity. Consider the linear elastic equation $\operatorname{div}(\mathbb{C}e(u))=0$ in a domain $\Omega\subset \mathbb{R}^3$ with smooth boundary, where Dirichlet and Neumann boundary conditions are separated by a finite number of smooth, nonintersecting closed curves. Here, $\mathbb{C}$ denotes the elasticity tensor and $e(u)$ the symmetrized gradient. Then $ u\in W^{1,4-\varepsilon}(\Omega)$ for any $\varepsilon>0$. For details, we refer to Theorem 8.1.7 and §8.3.1 in \cite{MaRO2010EEPD}. Lastly, the given bounds on $|\Re \lambda|$ are sharp in the sense that there exist a sequence of elliptic systems for which the corresponding $\lambda$ approaches these bounds (Section \ref{herewediscuss}).\medskip

While the approach provides new insights and results, its limitations must be acknowledged. This work focuses exclusively on the model problem in a planar angle. Hence, the results are applicable to two-dimensional domains and three-dimensional domains $\Omega\subset \mathbb{R}^3$, where a change of boundary conditions on $\partial \Omega$ manifests as an edge. Other geometric structures, such as vertices, are not covered and will be part of a subsequent paper.

\subsection*{Structure of the paper}
In \Cref{the model problem}, following \cite{MaRO2010EEPD} and \cite{VMR2001SPCS}, the model problem is introduced for an elliptic operator of the form
\begin{align*}
    L_A(\partial_x,\partial_y)=A_{11}\partial_x^2+2 A_{12}\partial_x\partial_y +A_{22}\partial_y^2
\end{align*}
for $A_\bullet\in \Mat_\ell(\mathbb{R})$ and $\ell\in \mathbb{N}$ the dimension of the system. It is assumed that $L_A$ is \textit{strongly elliptic} (weaker than the \textit{formal positivity} condition in elasticity, see Lemma \ref{seemeti}). The domain is given by the two-dimensional angle $\mathcal{K}_\alpha$ (\ref{kaha}) for some $0<\alpha\leq 2\pi$. We consider either Dirichlet or Neumann boundary conditions on the two sides of the angle, respectively. Using a factorization result for matrix polynomials in \cite{GLRMaPo2009}, Lemma \ref{posdef} yields a decomposition
\begin{align*}
    L_A(\xi_1,\xi_2)=A_{22}^{1/2}(V^* \xi_1 - \Id_\ell \xi_2)(V\xi_1-\Id_\ell \xi_2) A_{22}^{1/2}\quad \forall (\xi_1,\xi_2)\in \mathbb{C}^2
\end{align*}
for $V$ a complex-valued matrix with $\sigma(V)\subset \{z\in \mathbb{C}:\Im z>0\}$. The matrix $V$ is unique (Theorem \ref{algtheo}), and is referred to as the \textit{standard root} of $L_A$.
\smallskip

In \Cref{analysis of the model problem without boundary conditions}, it is shown that it suffices to consider elliptic systems with $A_{22}=\Id$, referred to as \textit{monic elliptic systems}. Using $V$, one can show that all solutions of the model problem without boundary conditions are given by (Prop. \ref{prar}):
\begin{align}\label{forthesol}
    u_{\lambda}:\mathbb{R}^2\setminus \{0\}\to \mathbb{C}^\ell,~ (x_1,x_2)\mapsto (x_1 \Id_\ell+x_2 V)^{\lambda}c_1+(x_1 \Id_\ell +x_2\overline{V})^{\lambda} c_2,\quad c_1,c_2\in \mathbb{C}^\ell.
\end{align}
The exponentiation of matrices here is defined via the \textit{functional calculus}, and the choice of complex exponentiation $\bullet\mapsto \bullet^\lambda$ required for smooth solutions is discussed. The constants $c_1$, $c_2$ are to be determined by the boundary conditions.
\smallskip

In \Cref{analysis of the model problem}, Dirichlet, mixed, and Neumann boundary conditions for $u_\lambda$ are implemented. It is shown that existence of a nontrivial solution $r^\lambda v$ for the model problem with angle $\alpha$ is for each boundary condition equivalent to $0\in \sigma(M_{\lambda,\alpha})$ for some $M_{\lambda,\alpha}\in \Mat_\ell(\mathbb{C})$ (Prop. \ref{dirstate}, Prop. \ref{mixstate}, Prop. \ref{neumannstate}). E.g., for Dirichlet boundary conditions, we get the condition
\begin{align}\label{zusammenhang herstellen}
    0\in \sigma(M_{\lambda,\alpha})\quad \text{for }M_{\lambda,\alpha}= Z_\alpha^\lambda-\overline{Z_\alpha}^\lambda,
\end{align}
where $Z_\alpha\in \Mat_\ell(\mathbb{C})$ is a symmetric matrix derived from $V$. The matrices $M_{\lambda, \alpha}$ for mixed and Neumann boundary conditions have a similar structure.  Also, we introduce two ellipticity conditions, \textit{Neumann well-posedness} and \textit{contractive Neumann well-posedness} (Def. \ref{neumann and contr neumann}), related to spectral properties of $V$. The former is equivalent to the \textit{complementing boundary condition} for Agmon-Douglis-Nirenberg (ADN)-elliptic systems \cite{ADN1964EBSE} implementing Neumann boundary conditions (see Appendix \ref{ellli}). The latter is crucial for the statement of the main result (Theorem \ref{main mixed}) and relates to path-connectedness of Neumann well-posed systems to the Laplacian (Lemma \ref{pathconnext}).\smallskip

In \Cref{matrix eq associated}, utilizing the numerical range and results on fractional powers of accretive operators \cite{Haas2003FCSO}, we are able to provide bounds on the spectrum of matrices $M_{\lambda}\in \Mat_\ell(\mathbb{C})$ reflecting the structure of $M_{\lambda,\alpha}$ for Dirichlet (Theorem \ref{dir}) and mixed boundary conditions (Theorems \ref{mix}, \ref{mix2}).\smallskip

In \Cref{regularity results of the model problem}, using the results of the previous sections, we prove bounds on $|\Re \lambda|$ for solutions $r^\lambda v$ of the model problem subject to Dirichlet and mixed boundary conditions. In particular, the main result, Theorem \ref{main mixed}, establishes for mixed boundary conditions the bounds $ |\Re \lambda|\geq \frac{1}{2}$ for $\alpha\leq \pi$ and $| \Re \lambda|\geq \frac{1}{4}$ for $\alpha\leq 2\pi$, provided that the system is contractive Neumann well-posed. If the system is not contractive Neumann well-posed, then cases with $|\Re \lambda|<\frac{1}{2}$, resp. $|\Re \lambda|<\frac{1}{4}$, may occur, but only in the form $\Re \lambda=0$. In \Cref{summary}, Neumann boundary conditions, optimality of the given bounds, and the scalar case $\ell=1$ are briefly discussed.
\smallskip

In Appendix \ref{ellli}, the ellipticity conditions (contractive) Neumann well-posedness, formal positivity, and the complementing boundary condition are related to each other. In Appendix \ref{spundfunc}, the functional calculus is summarized. In Appendix \ref{numerrange and accreive}, results on accretive operators are adapted to our setting.

\begin{figure}[h!]
    \centering
    \includegraphics[width=0.6\linewidth]{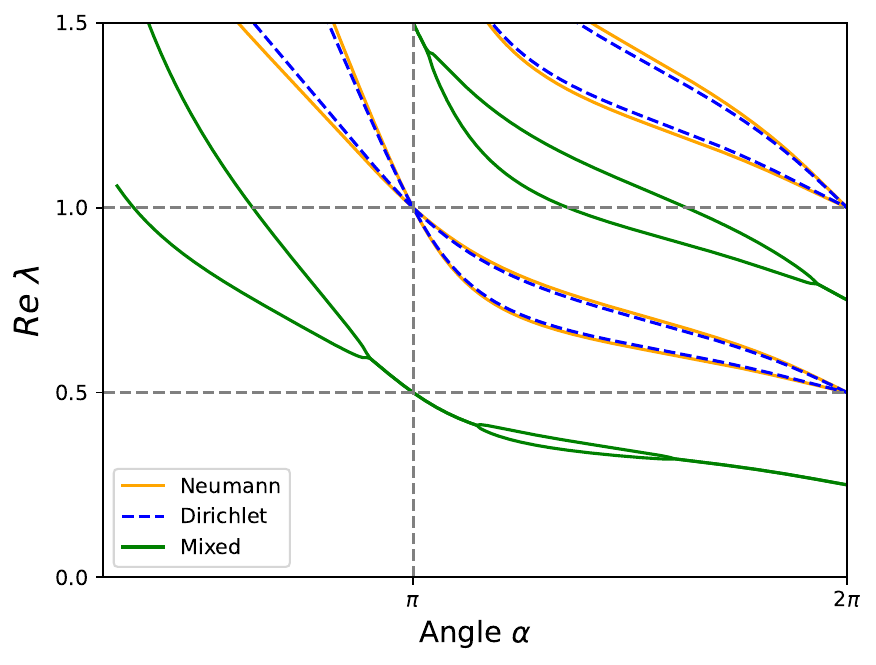}
    \caption{\small Relation between $\Re \lambda$ and $\alpha\in [1,2\pi]$ for different boundary conditions. The elliptic tuple is defined by $A_{11}=\begin{pmatrix}
        5&0.6\\0.6&1.5
    \end{pmatrix}$, $A_{12}=\begin{pmatrix}
        0.25&-0.4\\-0.4&-0.2
    \end{pmatrix}$, $A_{22}=\begin{pmatrix}
        1&0\\0&1
    \end{pmatrix}$. The branches for Dirichlet and Neumann boundary conditions are very close to each other.}
    \label{fig:1}
\end{figure}

\section{Notation and prerequisites}
LHS and RHS will be used as abbreviations for "left-hand side" and "right-hand side", respectively. For $r \in \mathbb{R} \setminus \{0\}$, let $\operatorname{sgn}(r) \in \{-1, +1\}$ denote the sign of $r$. $\mathbb{R}_{>0}$ denotes positive numbers. We use $(x_1,x_2)\in \mathbb{R}^2$ for Cartesian coordinates and $(r,\varphi)\in \mathbb{R}_{>0}\times [0,2\pi)$ for polar coordinates. $\Re z$ and $\Im z$ denote real and imaginary part of a complex number $z\in \mathbb{C}$ (or matrix), and $\overline{z}$ complex conjugation. For a set $A\subset \mathbb{C}$, we write $\operatorname{clos}(A)$ for its closure. Let us denote the (open) upper half-plane $\operatorname{UHP} \subset \mathbb{C}$ and the right half-plane $\operatorname{RHP}\subset \mathbb{C}$ by
\begin{align*}
    \operatorname{UHP}:=\{z\in \mathbb{C}: \Im(z)>0\},\quad \operatorname{RHP}:=\{z\in \mathbb{C}: \Re(z)>0\}.
\end{align*}
The lower half-plane and left half-plane are simply denoted by $-\operatorname{UHP}$ and $-\operatorname{RHP}$, where we use the following notation for set operations: For two sets $A,B\subset \mathbb{C}$ and $r\in \mathbb{C}$ we write
\begin{align*}
    A+B&:=\{a+b:a\in A,~b\in B\},\quad rA:=\{r\cdot a: a\in A\},\quad \overline{A}:=\{\overline{a}:a\in A\},\\ 
    \operatorname{Re}A&:=\{\Re(a):a\in A\},\quad  \operatorname{Im}A:=\{\Im(a):a\in A\}.
\end{align*}

Let us fix $\ell\in \mathbb{N}$. By $\langle \bullet,\bullet \rangle$, we denote the scalar product on $\mathbb{C}^\ell$, and by $\|\bullet\|$ the vector norm $\|v\|=\sqrt{\langle v,v\rangle}$ for $v\in \mathbb{C}^\ell$. Write $\Mat_\ell(\mathbb{C})$ for matrices of size $\ell\times \ell$ with entries in $\mathbb{C}$ and $ \Mat_\ell(\mathbb{R})$ if the entries are real-valued. $\Id_\ell\in \Mat_\ell(\mathbb{C})$ denotes the identity matrix. For $A\in \Mat_\ell(\mathbb{C})$, we write $A^T$ for the transpose, $A^*=\overline{A}^T$ for the adjoint, and $A^{-1}$ for the inverse (if it exists). We consider symmetric matrices $A^T=A$, Hermitian matrices $A^*=A$, and unitary matrices $A^*=A^{-1}$. The operator norm on $\Mat_\ell(\mathbb{C})$ is given by $\|A\| = \sup_{\|v\| = 1} \|Av\|$. The commutator is denoted by $[A,B]=AB-BA$. The set of all eigenvalues, the spectrum of $A\in \Mat_\ell(\mathbb{C})$, is denoted $\sigma(A)$. The spectral radius of $A$ is given by $\rho(A)=\sup_{\lambda \in \sigma(A)}|\lambda|$.
It satisfies $\rho(A)\leq \|A\|$, and for Hermitian matrices, $\rho(A)=\|A\|$. We call $A$ positive definite, and write $A>0$, if it is Hermitian and satisfies $\langle v,A v\rangle >0$ for any vector $v\in \mathbb{C}^\ell \setminus \{0\}$. A Hermitian matrix $A$ is positive definite if and only if $\sigma(A)\subset \mathbb{R}_{>0}$. If $\sigma(A)\subset \mathbb{R}_{> 0}\cup \{0\}$, $A$ is called positive semi-definite. The product $A^*A$ is always positive semi-definite, and it is positive definite if $0\notin \sigma(A)$. Any positive definite matrix $A$ has a unique positive definite square root denoted by $A^{1/2}$. For details, the reader is referred to \cite{Horn1985MaAn}. 

\section{The model problem}\label{the model problem}
This section introduces the model problem as presented in §6 of \cite{MaRO2010EEPD}. The domain of interest is given, for $0<\alpha\leq 2\pi$, by the two-dimensional angle
\begin{align}\label{kaha}
    \mathcal{K}_{\alpha}:=\{(r\cos(\varphi),r\sin(\varphi)): r>0,~ 0\leq \varphi \leq \alpha\}\subset \mathbb{R}^2,
\end{align}
which has the boundary $\partial \mathcal{K}_{\alpha}=\Gamma^-\cup \Gamma^+\cup \{(0,0)\}$ for
\begin{align*}
    \Gamma^-:=~\{(r,0):r> 0\}\quad \text{and}\quad \Gamma^+:=\{(r\cos(\alpha),r\sin(\alpha)):r>  0\}.
\end{align*}

For $A=(A_{11},A_{12},A_{22})$, where $A_\bullet\in \Mat_\ell(\mathbb{R})$ are symmetric matrices and $\ell\in \mathbb{N}$ the size of the system, the second-order differential operator $L_A$ is given by
\begin{align}\label{oftheform}
    L_A(\partial_{x_1},\partial_{x_2})&:=\sum_{i,j=1}^2 A_{ij}\partial_{x_i}\partial_{x_j}=A_{11}\partial^2_{x_1}+2A_{12}\partial_{x_1}\partial_{x_2}+A_{22}\partial_{x_2}^2,
\end{align}
where we set $A_{21}=A_{12}$ in the following. The conormal derivatives
$N_A^\pm$ associated to $L_A$ on $\Gamma^\pm$ are given by
\begin{align*}
    N^-_{A}(\partial_{x_1},\partial_{x_2})=& ~ N_{A}(0,\partial_{x_1},\partial_{x_2})\quad \text{and}\quad 
    N^+_{A}(\partial_{x_1},\partial_{x_2})=  N_{A}(\alpha,\partial_{x_1},\partial_{x_2})\quad \text{for}\\
     N_{A}(\varphi,\partial_{x_1},\partial_{x_2}):=&\sum_{i,j=1}^2 A_{ij}n_i\partial_{x_j}=A_{11} n_1\partial_{x_1}+A_{12}(n_1\partial_{x_2}+n_2\partial_{x_1})+A_{22}n_2 \partial_{x_2},\quad \varphi \in \{0,\alpha \},
\end{align*}
where $n=(-\sin(\varphi),\cos(\varphi))$ for $\varphi \in \{0,\alpha \}$ is the normal vector perpendicular to $\Gamma^\pm$. We investigate vector-valued solutions $u:\mathcal{K}_{\alpha}\to \mathbb{C}^\ell$ to the equations
\begin{align}\label{soe}
    L_A(\partial_{x_1},\partial_{x_2}) u=0\quad \text{on } \mathcal{K}_{\alpha},\quad
    B_A^\pm(\partial_{x_1},\partial_{x_2}) u=0 \quad \text{on }\Gamma^\pm,
\end{align}
where $u$ can be decomposed in the \textit{radial form} $u(r,\varphi)=r^\lambda v(\varphi)$ for some $\lambda \in \mathbb{C}$ and $v:[0,\alpha]\to \mathbb{C}^\ell$ smooth. Here,
\begin{align}\label{ofthenorm}
    B_A^\pm(\partial_{x_1},\partial_{x_2})u:=(1-d^\pm)u+d^\pm N^\pm_{A}(\partial_{x_1},\partial_{x_2})u
\end{align}
for $d^\pm\in \{0,1\}$ such that $(d^+,d^-)=(0,0)$ implements Dirichlet, $(d^+,d^-)=(1,1)$ Neumann, and $(d^+,d^-)=(0,1)$ mixed boundary conditions. The problem of finding a solution $u$ of radial form to (\ref{soe}) is called the \textit{model problem}. If we only consider
\begin{align}\label{soeo}
    L_A(\partial_{x_1},\partial_{x_2}) u=0\quad \text{on } \mathcal{K}_{\alpha},
\end{align}
then we call (\ref{soeo}) the \textit{model problem without boundary conditions}. Our question is the following: For fixed $\alpha\in (0,2\pi]$, for which $\lambda\in \mathbb{C}$ can we expect a solution of the form $r^\lambda v$ to the model problem (\ref{soe})? More specifically, what is the smallest value for $|\Re \lambda|$ that we can expect for a solution? These results can be translated to regularity of solutions for strongly elliptic systems in polyhedral domains (see §2 and §6 in \cite{MaRO2010EEPD}).

\subsection{Ellipticity}

We use the following notions of ellipticity for $L_A$ which can be found in §1.1.2 of \cite{MaRO2010EEPD}.
 \begin{definition}\label{the conda}
     Consider the differential operator $L_A$ in (\ref{oftheform}) for $A_\bullet\in \Mat_\ell(\mathbb{R})$ symmetric. For $\xi=(\xi_1,\xi_2)\in \mathbb{R}^2$, we define the polynomial $L_A(\xi):=\sum_{i,j=1}^2 A_{ij}~\xi_i\xi_j$. We call $L_A$ \textit{elliptic} if
         \begin{align}\label{nostra}
             \det(L_A(\xi))\neq 0\quad \forall \xi \in \mathbb{R}^2 \setminus \{0\},
         \end{align}
and say $L_A$ is \textit{strongly elliptic} if there exists $\kappa>0$ such that
\begin{align}\label{kal}
   \langle L_A(\xi) \eta,\eta\rangle \geq \kappa \|\eta\|^2 \|\xi\|^2\quad \forall \eta\in \mathbb{C}^\ell,~\xi\in \mathbb{R}^2.
\end{align}
 \end{definition}

\begin{remark}
   i) $L_A(\xi)\in \Mat_\ell(\mathbb{R})$ is symmetric real-valued and thus Hermitian for any $\xi\in \mathbb{R}^2$. So the LHS in (\ref{kal}) is always real.\\
   ii) The condition for strong ellipticity is sometimes referred to as \textit{Legendre-Hadamard condition}.
\end{remark}

Obviously, strong ellipticity implies ellipticity. From now on, we always assume that $A_{11}$ and $A_{22}$ are positive definite (which is necessary, but not sufficient for strong ellipticity). In this case, one can show that the notions coincide, simplifying the analysis.
    
\begin{lemma}\label{posdef}
    Assume that $A_{11},A_{12},A_{22}\in \Mat_\ell(\mathbb{R})$ are symmetric matrices and $A_{11}$, $A_{22}$ are positive definite. Then the following are equivalent:
    \begin{enumerate}
        \item[i)] $\det(L_A(\xi))\neq 0$ for all $\xi\in \mathbb{R}^2$ of the form $\xi=(1,\beta)\in \mathbb{R}^2$.
        \item[ii)] $L_A$ is elliptic.
        \item[iii)] $L_A$ is strongly elliptic.
    \end{enumerate}
    In these cases, there exists the factorization
    \begin{align}\label{factobacto}
        L_A(\xi_1,\xi_2)=A_{22}^{1/2}(V^* \xi_1 - \Id_\ell \xi_2)(V\xi_1-\Id_\ell \xi_2) A_{22}^{1/2}\quad \forall (\xi_1,\xi_2)\in \mathbb{C}^2
    \end{align}
    for some $V\in \Mat_\ell(\mathbb{C})$ with $\sigma(V)\subset \operatorname{UHP} $.
\end{lemma}

\begin{proof}
The implications iii)$\implies$ ii) $\implies$ i) are clear. Now, assume $\det( L_A(1,\beta))\neq 0$ for all $\beta \in \mathbb{R}$ and show that $L_A$ is strongly elliptic.

\smallskip
    
First, we show $L_A(1,\beta)>0$ for $\beta\in \mathbb{R}$. Symmetry of $L_A(1,\beta)$ follows from symmetry of the $A_\bullet$'s. Since $\det(L_A(1,\beta))$ is the product of eigenvalues of $L_A(1,\beta)$, and $L_A(1,0)=A_{11}>0$, continuity of eigenvalues ensures  $\sigma(L_A(1,\beta))\subset \mathbb{R}_{>0}$ for all $\beta\in \mathbb{R}$.
\smallskip

Note that one can write for any $\beta\in \mathbb{C}$:
    \begin{align}\label{monica}
        L_A(1,\beta)=A_{22}^{1/2}L_{\tilde{A}}(1,\beta)A^{1/2}_{22},\quad \text{with}\quad \tilde{A}:=(A_{22}^{-1/2}A_{11} A_{22}^{-1/2},A_{22}^{-1/2}A_{12} A_{22}^{-1/2},\Id_\ell).
    \end{align}

Sylvester's law of inertia implies $L_{\tilde{A}}(1,\beta)>0$ for $\beta\in \mathbb{R}$. Due to Thm. 12.8 in \cite{GLRMaPo2009}, there is $V\in \Mat_\ell(\mathbb{C})$ with $\sigma(V)\subset \operatorname{clos}(\operatorname{UHP})$ such that 
\begin{align}\label{märzstatt2}
    L_{\tilde{A}}(1,\beta)=(V^* -\Id_\ell \beta)(V-\Id_\ell \beta)\quad \forall \beta\in \mathbb{C}.
\end{align}
In fact, $\sigma(V)\subset \operatorname{UHP}$ because if $\beta \in \sigma(V)\cap  \mathbb{R}$, then $\det( L_{\tilde{A}}(1,\beta))= 0=\det( L_A(1,\beta))$, contradicting ellipticity. The factorization (\ref{factobacto}) is deduced from (\ref{monica}) and (\ref{märzstatt2}).
\smallskip

Next, define the compact set $K:=\{(\xi,\eta)\in \mathbb{R}^2\times \mathbb{C}^\ell : \|\xi\|=1=\|\eta\|\}$. Note that
\begin{align*}
    \langle L_A(\xi)\eta,\eta\rangle=\| (\xi_1 V-\xi_2 \Id_\ell)A_{22}^{1/2}\eta\|^2
\end{align*}
and set $\kappa:=\inf_{K} \| (\xi_1 V-\xi_2 \Id_\ell)A_{22}^{1/2}\eta\|^2$. Due to compactness of $K$ and $\sigma(V)\cap \mathbb{R}\neq \emptyset$, $\kappa>0$, which establishes (\ref{kal}).
\end{proof}

Throughout this work, we refer to $A$ as an \textit{elliptic system} or \textit{elliptic tuple} if $A=(A_{11},A_{12},A_{22})$ are all symmetric, $A_{11}>0$, $A_{22}>0$, and $L_A$ is elliptic. We refer to elliptic tuples $A$ with $A_{22}=\Id_\ell$ as \textit{monic elliptic tuples} and to $L_A$ as a \textit{monic elliptic operator}. For an elliptic tuple $A$ we call the tuple $\tilde{A}$ in (\ref{monica}) its \textit{monic reduction}. By i) in Lemma \ref{posdef}, it is clear that $L_{\tilde{A}}$ is elliptic. Also, for $L_A$ a monic elliptic operator, we call $V\in\Mat_\ell(\mathbb{C})$ fulfilling (\ref{märzstatt2}) and $\sigma(V)\subset \operatorname{UHP}$ a \textit{standard root} of $L_A$. By its definition, the standard root of $L_A$ is the same as for $L_{\tilde{A}}$. The investigation of standard roots will be crucial to derive solutions to the model problem.

\begin{example}\label{C=0,laplace}
    For the Laplacian tuple $A=(\Id_\ell,0,\Id_\ell)$, a standard root is given by $V=i\Id_\ell$.
\end{example}

Thus far, we have only described the ellipticity of the operator $L_A$. However, to ensure the well-posedness of the elliptic problem, it is also necessary to specify a \textit{complementing condition} for the boundary operators $B_A^\pm$. These conditions, introduced in the framework of elliptic systems by Agmon, Douglis, and Nirenberg in \cite{ADN1964EBSE}, establish compatibility between the boundary operators and the elliptic operator. The compatibility condition will naturally arise in our analysis (Section \ref{neumann boundary yea}), and we postpone a detailed discussion until that point.

\subsection{The operator pencil}\label{eigenvalofop}
In this section, the model problem is translated to a parameter-dependent second-order ODE, and the operator pencil is introduced. The reference is §6.1.3 of \cite{MaRO2010EEPD}. Writing $u=r^\lambda v$ for $\lambda\in \mathbb{C}$, we define the $\lambda$-dependent differential operators $\mathcal{L}_A(\partial_\varphi,\lambda)$, $\mathcal{B}_A^\pm(\partial_\varphi,\lambda)$ and $\mathcal{N}_A(\partial_\varphi,\lambda)$ by:
\begin{align*}
    \mathcal{L}_A(\partial_\varphi,\lambda) v &:=r^{2-\lambda}L_A(\partial)r^\lambda v,\\
    \mathcal{N}_A(\partial_\varphi,\lambda) v&:=r^{1-\lambda}N_{A}(\partial)r^\lambda v,\\
     \mathcal{B}_A^\pm(\partial_\varphi,\lambda)v&:=(1-d^\pm)v+d^\pm \mathcal{N}_A(\partial_\varphi,\lambda) v.
\end{align*}
A long but straightforward calculation (compare also to Def. 6 in \cite{HDKREMPM2008}) shows that:
\begin{align}\label{L__A}
     \mathcal{L}_A(\partial_\varphi,\lambda)&=b_2(\varphi)\partial_\varphi^2+(\lambda-1) b_1(\varphi)\partial_\varphi+\lambda(\lambda-1) b_0(\varphi)+\lambda b_2(\varphi),  \\ \nonumber
    \mathcal{N}_A(\partial_\varphi,\lambda)&=b_2(\varphi)\partial_\varphi +\frac{\lambda}{2} b_1(\varphi),
\end{align}
where the $b_\bullet$'s are the periodic functions:
\begin{align}\label{hjh}
      b_0(\varphi)&=A_{11}\cos(\varphi)^2+A_{22}\sin(\varphi)^2 +2A_{12}\sin(\varphi)\cos(\varphi),\\ \nonumber
      b_1(\varphi)&=2(A_{22}-A_{11})\sin(\varphi)\cos(\varphi)+2A_{12}(\cos(\varphi)^2-\sin(\varphi)^2),\\ \nonumber
    b_2(\varphi)&=A_{11}\sin(\varphi)^2+A_{22}\cos(\varphi)^2-2A_{12}\cos(\varphi)\sin(\varphi).
\end{align} 
We define the $\lambda$-dependent mapping:
\begin{align}\label{alamb}
    \mathcal{A}(\lambda):W^{2,2}((0,\alpha),\mathbb{C}^\ell)\to L^2((0,\alpha),\mathbb{C}^\ell)\times \mathbb{C}^\ell\times \mathbb{C}^\ell
\end{align}
by
\begin{equation*}
    v \mapsto \left( \mathcal{L}_A(\lambda) v,~ \mathcal{B}_A^-(\lambda) v \big|_{\varphi = 0}, ~\mathcal{B}_A^+(\lambda) v \big|_{\varphi = \alpha} \right).
\end{equation*}
Here, $W^{2,2}$ and $L^2$ denote the usual Sobolev and Lebesgue space. In the literature, $\mathcal{A}(\lambda)$ is called the \textit{operator pencil}. If there exist $\lambda\in \mathbb{C}$ and $v\neq 0$ such that $\mathcal{A}(\lambda)v=0$, then $\lambda$ is called an \textit{eigenvalue} of $\mathcal{A}$ and $v$ an \textit{eigenvector} to $\lambda$. See  §1 in \cite{VMR2001SPCS} for an introduction to operator pencils and further references. With the above derivations, the model problem is reduced to a $\lambda$-dependent second-order ODE, such that: \textit{The model problem (\ref{soe}) has a solution of the form $r^\lambda v$ for $\lambda \in \mathbb{C}$ if and only if $\lambda$ is an eigenvalue with eigenvector $v$ of the operator pencil $\mathcal{A}(\lambda)$ in (\ref{alamb}).}

\begin{remark}[Fundamental solutions for $\mathcal{L}_A(\partial_\varphi,\lambda)v=0$]
Fix $\lambda\in \mathbb{C}$, and observe that the system of $\ell$ second-order ODEs given by $\mathcal{L}_A(\partial_\varphi,\lambda)v=0$ can be reduced, by a standard trick, to a system of $2\ell$ first-order ODEs of the form $\partial_\varphi y=M(\varphi,\lambda) y$. Here, the entries of $M(\varphi,\lambda)\in \Mat_{2\ell}(\mathbb{C})$ are analytic in $\varphi$. The reduction requires inverting $b_2(\varphi)$, which is possible since $b_2(\varphi)=L_A(\xi)>0$ for $\xi=(\sin(\varphi),-\cos(\varphi))\in \mathbb{R}^2$ by Lemma \ref{posdef}. Using results from §IV.10 in \cite{Hart2002ODE}, the system of first-order ODEs has a fundamental matrix $Y(\varphi,\lambda)$ which is analytic in $\varphi$. This implies that the set of solutions to $\mathcal{L}_A(\partial_\varphi,\lambda)v=0$ spans a $2\ell$-dimensional complex vectorspace. By the preceding discussion, the same holds true for the set of solutions $u_\lambda=r^\lambda v$ to the model problem without boundary conditions. The analytic dependence of $v(\lambda)$ on $\lambda\in \mathbb{C}$ is revealed later (Section \ref{explisolf}).
\end{remark}

\begin{example}[Laplace equation]
    If we assume $A_{11}=\Id_\ell=A_{22}$ and $A_{12}=0$, then the system reduces to decoupled Laplace equations in $\ell$ components. In this case, $\mathcal{L}_A(\lambda)=\Id_\ell (\partial_\varphi^2+\lambda^2)$ and solutions to the $\lambda$-dependent ODE without boundary conditions are given by
\begin{align*}
    v_\lambda(\varphi)=c_1 \sin(\lambda \varphi)+c_2 \cos(\lambda \varphi)\quad \text{for }c_1,c_2\in \mathbb{C}^\ell.
\end{align*}
Let us implement boundary conditions for $\varphi\in \{0,\alpha\}$. Note that $\mathcal{N}_A(\lambda)=\Id_\ell \partial_\varphi$, so nontrivial solutions are given by $c_2=0$ and $\lambda\in \frac{\pi}{\alpha}\cdot \mathbb{Z}\setminus \{0\}$ for Dirichlet, $c_1=0$ and $\lambda\in \frac{\pi}{\alpha}\cdot \mathbb{Z}$  for Neumann, and $c_2=0$ and $\lambda\in \frac{\pi}{2\alpha} \cdot \mathbb{Z}\setminus \{0\}$ for mixed boundary conditions (compare to §2.1 in \cite{VMR2001SPCS}). This leads for $\alpha \leq \pi$ to the bounds $|\Re \lambda|\geq 1$ for Dirichlet and Neumann boundary conditions (ignoring constant solutions at $\lambda=0$ for the latter) and $|\Re \lambda |\geq \frac{1}{2}$ for mixed boundary conditions. Although this is the simplest elliptic system, we will derive similar lower bounds for more general systems.
\end{example}

\subsection{The case $\lambda=0$}
The subsequent derivation does not cover $\lambda= 0$, which is why we address this case now. A solution $u=r^\lambda v$ to the model problem (\ref{soe}) is called a \textit{trivial solution} if $v=0$. The following is not a new result but included for completeness.
\begin{lemma}\label{lambda=0}
    The model problem (\ref{soe}) admits for $\lambda=0$ and any angle $0<\alpha\leq 2\pi$ only the trivial solution for Dirichlet and mixed boundary conditions, and only constant solutions for Neumann boundary conditions.
\end{lemma}
\begin{proof}
Consider an elliptic tuple $A$, and assume $\mathcal{L}_A(\lambda)v=0$ for $\lambda=0$. By (\ref{L__A}), this reduces to
    \begin{align}\label{generaleeass}
    \partial_\varphi(b_2 \partial_\varphi v)(\varphi) =0,
    \end{align}
    due to $\partial_\varphi b_2=-b_1$. The $2\ell$ linearly independent solutions are given by
    \begin{align*}
        v(\varphi)=c_1+ P(\varphi)c_2 \quad \text{for }c_\bullet \in \mathbb{C}^\ell,
    \end{align*}
where $P(\varphi)=\int^\varphi_0 b_2^{-1}(s)ds $. Note that $P(\varphi)>0$ for $\varphi>0$, ensuring that $P(\varphi)$ is invertible. Dirichlet boundary conditions yield 
\begin{align*}
    c_1=0\quad \text{and}\quad c_1+P(\alpha)c_2=0\implies c_1=0=c_2,
\end{align*}
so only the trivial solution exists. Neumann boundary conditions (check (\ref{L__A})) yield  $c_2=0$ for $\varphi \in \{0,\alpha\}$, so $v(\varphi)=c_1$ is the most general solution. Finally, for mixed boundary conditions, both $c_1=0=c_2$ are enforced, so only the trival solution exists.
\end{proof}

\section{Analysis of the model problem without boundary conditions}\label{analysis of the model problem without boundary conditions}
\subsection{Algebraic properties of $L_A$}
The next result reduces the discussion to monic elliptic tuples.

\begin{lemma}\label{reducioo}
Consider an elliptic tuple $A=(A_{11},A_{12},A_{22})$ and its monic reduction $\tilde{A}$. The model problem (\ref{soe}) admits for $(L_A,B_A)$ the solution $u=r^\lambda v$ if and only if the model problem for $(L_{\tilde{A}},B_{\tilde{A}})$ admits the solution $\tilde{u}=r^\lambda \tilde{v}$, where $\tilde{v}=A_{22}^{1/2}v$.
\end{lemma}

\begin{proof}
Assume that  $L_A u=0$. Clearly, $\tilde{u}=A^{1/2}_{22} u$ is a solution to $L_{\tilde{A}} \tilde{u}=0$ by (\ref{monica}). Also, if $u=r^\lambda v$, then $\tilde{u}=r^\lambda \tilde{v}$ for $\tilde{v}=A^{1/2}_{22} v$. For Dirichlet boundary conditions, we have $v(\varphi)=0$ if and only if $A_{22}^{1/2}v(\varphi)=0$. For Neumann boundary conditions, one uses $\mathcal{N}_A(\partial_\varphi,\lambda) v (\varphi)= A_{22}^{1/2} \mathcal{N}_{\tilde{A}}(\partial_\varphi,\lambda) A^{1/2}_{22}v(\varphi) $. Since $A_{22}^{1/2}$ is invertible, the converse implication is clear.
\end{proof}
For the remainder, we reduce the discussion to monic elliptic tuples $A=(A_{11},A_{12},\Id_\ell)$. Recall that, by Lemma \ref{posdef}, there exists a standard zero $V\in \Mat_\ell(\mathbb{C})$ with $\sigma(V)\subset \operatorname{UHP}$ such that
\begin{align}\label{jkj}
        L_A(1,\beta)=A_{11}+2\beta A_{12}+\beta^2 \Id_\ell=(V^* -\Id_\ell \beta)(V-\Id_\ell \beta).
    \end{align}

\begin{lemma}\label{maybel}
1. Consider a monic elliptic tuple $A=(A_{11},A_{12},\Id_\ell)$. Assume $V=C+iD$ with $C=\Re V$ and $D=\Im V$ fulfills (\ref{jkj}). Then:
\begin{alignat}{3}\label{defino}
    A_{11}&=C^T C+D^TD,&\quad  A_{12}&=-\frac{1}{2}(C+C^T),\\ \label{defino2}
    D&=D^T,&\quad C^TD&=D^TC.
\end{alignat}
2. On the other hand, assume $V=C+iD$ for $C,D\in \Mat_\ell(\mathbb{R})$ which satisfy the algebraic relations (\ref{defino2}). Then:
\begin{enumerate}
    \item[i)] $V$ satisfies (\ref{jkj}) for $A_\bullet$ as defined by (\ref{defino}).
    \item[ii)] If $ \sigma(V)\cap \mathbb{R}=\emptyset $, then the tuple defined by (\ref{defino}) and $A_{22}=\Id_\ell$ is monic elliptic.
\end{enumerate}

\end{lemma}

\begin{proof}
1) Writing $V=C+iD$, $A_{11}=V^* V$, and $A_{12}=-\frac{1}{2}( V^* +V)$, we get:
\begin{align*}
     A_{11}=C^TC+D^TD+i(-D^TC+C^TD),\quad A_{12}=-\frac{1}{2}( C^T+C)+i(-D^T+D).
\end{align*}
The imaginary parts vanish since $A_{11},A_{12}\in \Mat_\ell(\mathbb{R})$. Thus, (\ref{defino}) and (\ref{defino2}) follow.\\
2.i) Assume $V=C+iD$ for $C,D\in \Mat_\ell(\mathbb{R})$ satisfies (\ref{defino2}). It is clear from the last calculation that we obtain (\ref{jkj}) for $A_\bullet$ given in (\ref{defino}).\\
2.ii) We additionally assume $ \sigma(V)\cap \mathbb{R}=\emptyset $ and prove that the corresponding operator $L_A$ is elliptic. First, we check symmetry of $A_{11},A_{12},\Id_\ell$ and $A_{11},\Id_\ell>0$. Symmetry is clear by the definition of $A_{11}$ and $A_{12}$ given in (\ref{defino}). For positive definiteness, note that $A_{11}=V^*V$, which is positive definite since $0\notin \sigma(V)$. Lastly, using Lemma \ref{posdef}, it suffices to show $\det(L_A(1,\beta))\neq 0$ for any $\beta \in \mathbb{R}$. This follows from 
\begin{align*}
    \det(L_A(1,\beta))=\det\big( (V^*-\Id_\ell  \beta)(V-\Id_\ell \beta)\big)=\det( V^*-\Id_\ell  \beta)\det (V-\Id_\ell \beta )\neq 0\quad \forall \beta\in \mathbb{R},
\end{align*}
since $\sigma(V)\cap \mathbb{R}=\emptyset $ and thus also $\sigma(V^*)\cap \mathbb{R}=\emptyset $.
\end{proof}
For a monic elliptic tuple $A$, there exist several matrices $V\in \Mat_\ell(\mathbb{C})$ that satisfy (\ref{jkj}). However, our goal is to show that $\sigma(V)\subset \operatorname{UHP}$ uniquely characterizes $V$, i.e., standard zeros are unique.

\begin{lemma}\label{inv}
    Let $V$ be a standard root of a monic elliptic tuple $A$. Then $D=\Im V$ is positive definite.
\end{lemma}
\begin{proof}
    Let us write $V=C+iD$ for $C,D\in \Mat_\ell(\mathbb{R})$, where $V$ is the standard root under consideration. Due to (\ref{defino2}), we have $D=D^T$ such that all eigenvalues of $D$ are real, and we need to show that they are all positive. First, we argue that $0\notin \sigma(D)$. Assume the contrary, $0\in \sigma(D)$, and derive a contradiction. In this case, there is $v\in \mathbb{R}^\ell \setminus \{0\}$ such that $Dv=0$. Due to the commutativity relations in (\ref{defino2}), we have $0=C^T Dv=DCv$. From this, it follows that either $Cv=0$ or $Cv\neq 0$ is also an eigenvector of $D$ with eigenvalue $0$. If $Cv=0$, it would follow that $V v=0$, which is a contradiction to $\sigma(V)\cap \mathbb{R}=\emptyset$. So, assuming the latter, $\operatorname{span}_\mathbb{\mathbb{C}}\left(\{v,Cv\}\right)\subset \ker D$. But repeating the discussion with $Cv$ instead of $v$, we conclude that $C^n v\neq 0$, for arbitrary $n\in \mathbb{N}$, must also be an eigenvector of $D$ with eigenvalue $0$. Define the cyclic (complex) subspace $ S_v:=\operatorname{span}_\mathbb{C}\left(\{C^n v:n\in \mathbb{N}_0\}\right)$ generated by $v$. By the above derivation, we have $D(S_v)=\{0\}$, and by its definition, it is clear that $C(S_v)\subset S_v$. This implies $C|_{S_v}=V|_{S_v}$. Thus, $S_v$ is an eigenspace of $V$, and moreover, $V$ has a real-valued matrix representation (the same as $C$) for the subspace $S_v$. This implies that $V$ has either a real eigenvalue or two different complex conjugated eigenvalues, which contradicts $\sigma(V)\subset \operatorname{UHP}$. Thus, we have shown that $0\notin \sigma(D)$.\\
Next, we argue why $D$ cannot have negative eigenvalues. This follows by a simple scaling argument. Define for any $\rho\geq 1$ the matrix $V_\rho:=C+i \cdot \rho D$. Observe that $C_\rho:=\Re V_\rho=C$ and $D_\rho:= \Im V_\rho=\rho D$ still fulfill the algebraic relations (\ref{defino2}). Moreover, the operator $L_{A,\rho}$ defined by $V_\rho$ (second part of Lemma \ref{maybel}) is elliptic for any $\rho\geq 1$, since we can write the matrix polynomial as
\begin{align*}
    L_{A,\rho}(1,\beta)=(\rho-1)D^2+A_{11}+2A_{12}\beta+\Id_\ell\beta^2=(\rho-1)D^2+L_A(1,\beta),
\end{align*}
and $L_{A}(1,\beta)> 0$, for $\beta\in \mathbb{R}$, and $(\rho-1)D^2$ is positive semi-definite. Now, let us assume that $r\in \sigma(D)$ for some $r<0$ and derive a contradiction. Since all eigenvalues of $D$ are nonzero, we get for the spectrum
    \begin{align*}
        \sigma(V_{\rho})=\sigma(C+i\rho D)\xrightarrow[]{\rho \to \infty} i\rho  \cdot \sigma(D)
    \end{align*}
    in the appropriate sense. Since $0>r\in \sigma(D)$, for sufficiently large $\rho>1$, there is $\beta_\rho \in \sigma(V_\rho)$ with $\Im \beta_\rho<0$. By continuity, and because $\sigma(V_1)\subset \operatorname{UHP}$, there must be an intermediate $1<\tilde{\rho}<\rho$ such that $\sigma(V_{\tilde{\rho}})\cap \mathbb{R}\neq \emptyset$, contradicting the ellipticity of $L_{A,\tilde{\rho}}$.
\end{proof}
The next statement is the main result for standard roots.

\begin{theorem}\label{algtheo}
Consider $L_A$, a monic elliptic operator, and $V$ a standard root of $L_A$.
\begin{enumerate}
    \item[i)] The standard root of $L_A$ is unique.
    \item[ii)] $V$ satisfies $A_{11}+2A_{12}V+\Id_\ell V^2=0=A_{11}+2A_{12}\overline{V}+\Id_\ell \overline{V}^2$.
    \item[iii)] $V$ can be written as $V=(S+i\Id_\ell)D$, where both $S,D\in \Mat_\ell(\mathbb{R})$ are symmetric and $D>0$.
\end{enumerate}
Conversely, any $V=(S+i\Id_\ell)D$ with $S,D\in \Mat_\ell(\mathbb{R})$ symmetric and $D>0$ is the standard root of the monic elliptic operator $L_A\,$ for the tuple
    \begin{align}\label{melookme}
        A=(A_{11}=V^*V,A_{12}=-\frac{1}{2}(V+V^*),A_{22}=\Id_\ell).
    \end{align}
\end{theorem}

\begin{proof}
Consider $L_A$, a monic elliptic operator, and let $V\in \Mat_\ell(\mathbb{C})$ be a standard root.

i) Uniqueness of $V$ follows from $L_A(1,\beta)=(V^*-\Id_\ell \beta)(V-\Id_\ell \beta)$, $\sigma(V)\subset \operatorname{UHP}$, $\sigma(V^*)\subset -\operatorname{UHP}$, and uniqueness of \textit{monic $\Gamma$-spectral right divisors} for matrix polynomials. For the latter, we refer to §4.1 in \cite{GLRMaPo2009}, in particular Theorem 4.1 and the comment thereafter.

ii) The first equality follows by $A_{11}=V^*V$ and $A_{12}=-\frac{1}{2}(V+V^*)$. The second by complex conjugation and $A_\bullet \in \Mat_\ell(\mathbb{R})$.

iii) By Lemma \ref{maybel} and Lemma \ref{inv}, we can write $V=C+iD$ for $C,D\in \Mat_\ell(\mathbb{R})$, where $D>0$ and (\ref{defino2}) holds. In particular, we can invert $D$, and thus $V=(S+i\Id_\ell)D$ for $S=CD^{-1}$. The symmetry of $S$ follows from:
    \begin{align*}
        S^T=(D^{-1})^T C^T=D^{-1} C^T D D^{-1}\stackrel{(\ref{defino2})}{=}D^{-1} D C D^{-1}=C D^{-1}=S.
    \end{align*}
Lastly, consider $V=(S+i\Id_\ell)D$ with $S,D\in \Mat(\mathbb{R},\ell)$ symmetric and $D>0$. We show that $V$ is the standard root of $L_A$ for $A$ given in (\ref{melookme}). By Lemma \ref{maybel}, it suffices to verify the algebraic relations
    \begin{align*}
        D=D^T,\quad (SD)^TD=D (SD),
    \end{align*}
    which are trivially fulfilled by symmetry of $S$ and $D$, and to show $\sigma(V)\cap \mathbb{R}=\emptyset$. The latter statement follows, due to $D>0$, from Lemma \ref{symm} (proven later) and $\sigma(V)=\sigma(D^{1/2}VD^{-1/2})$.
\end{proof}

\subsection{Solutions for the model problem without boundary conditions}\label{explisolf}
It is time to close the gap between algebra and analysis and justify the time spent on exploring the algebraic structure of $L_A$. In this section, solutions for the model problem without boundary conditions are derived. The representation of solutions is new, but is based on ideas in §2.2 of \cite{CoDauFSAC2001}. For this, we first need to define the complex exponent $\lambda \in \mathbb{C}$ of a complex number $z \in \mathbb{C}\setminus \{0\}$. We do so by using three branches of the complex logarithm, i.e. three different $\arg$-functions. To distinguish between them, we introduce the set of symbols $a\in \{o,+,-\}$ and use these as decoration. Explicitly, we define
\begin{alignat}{2} \label{explsolss}
    z^{\lambda_a}:=&\exp(\lambda \log_a (z))\quad \text{for}\quad \log_a(z):=\log(|z|)+i\lambda \arg_a(z).
\end{alignat}
Here, $\log(r)$ for $r\in \mathbb{R}$ is simply the standard logarithm for positive real numbers. The $\arg_a$-functions are uniquely determined by requiring that $\log_a$ inverts $\exp$ on the domain $\mathbb{C}\setminus \{0\}$ and by the following conditions:
\begin{align*}
    \arg_+(z)\in [0,2\pi),\quad \arg_o(z)\in (-\pi,\pi],\quad \arg_-(z)\in (-2\pi,0]\quad \forall z\in \mathbb{C}\setminus \{0\}.
\end{align*}
With this choice, $\log_o$, denoted by $\log$ in the following, represents the principal logarithm, and it exhibits a discontinuity along the branch cut of the negative real axis. The introduction of $\log_+$ and $\log_-$ serves to provide continuous extensions of the logarithm that avoid the discontinuity at $\arg(z) = \pi$, shifting the branch cut instead to the positive real axis. Note that $\arg_o=\arg_+$ on $\operatorname{UHP}$ and $\arg_o=\arg_-$ on $-\operatorname{UHP}$.
\smallskip

Consider $0<\alpha<2\pi$ and the model problem for a monic elliptic tuple $A$ with standard root $V$. We define, for $\lambda \in \mathbb{C}\setminus \{0\}$ and $c_1,c_2\in \mathbb{C}^\ell$ arbitrary, the vector-valued functions:
\begin{align}\label{wunderdich}
    u_{\lambda}:\mathbb{R}^2\setminus \{0\}\to \mathbb{C}^\ell,~ (x_1,x_2)\mapsto (x_1 \Id_\ell+x_2 V)^{\lambda_+}c_1+(x_1 \Id_\ell +x_2\overline{V})^{\lambda_-} c_2.
\end{align}
Here, the exponentiation of matrices is defined via the functional calculus (see Appendix \ref{spundfunc}). It is well-defined if $0\notin \sigma(x_1\Id_\ell+x_2V)=x_1+x_2\sigma(V)$ is ensured for any $(x_1,x_2)\in \mathbb{R}^2\setminus \{0\}$, which follows from $\sigma(V)\cap \mathbb{R}=\emptyset$. We will show that any solution to the model problem without boundary conditions (\ref{soeo}) is of the form (\ref{wunderdich}). For this, note that $u_\lambda=r^\lambda v_\lambda$, where 
\begin{align}\label{wunderdich2}
     v_\lambda:[0,2\pi)\to \mathbb{C}^\ell,\quad \varphi \mapsto (\cos(\varphi) \Id_\ell+\sin(\varphi) V)^{\lambda_+}c_1+(\cos(\varphi) \Id_\ell +\sin(\varphi)\overline{V})^{\lambda_-} c_2.
\end{align}
\begin{remark}
    Let us briefly discuss the choice of $\lambda_\pm$. Since $\sigma(V)\subset \operatorname{UHP}$, we have
\begin{alignat*}{1}
    &\sigma(\cos(\varphi) \Id_\ell+\sin(\varphi) V) \subset\phantom{-} \operatorname{UHP} \quad\text{for } 0<\varphi<\pi,\\
    &\sigma(\cos(\varphi) \Id_\ell+\sin(\varphi) V) \subset -\operatorname{UHP} \quad\text{for } \pi<\varphi<2\pi.
\end{alignat*}
Because the $\arg$-function in the principle logarithm has a discontinuity at $\varphi=\pi$, the function $\varphi\mapsto(\cos(\varphi) \Id_\ell+\sin(\varphi) V)^{\lambda}$ is not continuous. This issue is resolved by changing $\lambda$ to $\lambda_+$. A similar reasoning applies to $\overline{V}$ and $\lambda_-$.
\end{remark}
Next, we show $L_A u_\lambda=0$. For this, observe that 
\begin{align*}
      &L_A(\partial_{x_1},\partial_{x_2})u_\lambda=\sum_{i,j=1}^2 A_{ij}\partial_{x_i}\partial_{x_j} u_\lambda\\
      =&\lambda(\lambda-1)\Big((A_{11}+2A_{12} V+V^2) (x_1 \Id_\ell+x_2 V)^{(\lambda-2)_+}c_1+(A_{11}+2A_{12} \overline{V}+\overline{V}^2) (x_1 \Id_\ell+x_2 \overline{V})^{(\lambda-2)_-}c_2\Big)
\end{align*}
vanishes by Theorem \ref{algtheo}. Here, the chain rule is applied, and it is important to note that the differentiation rules for $\bullet^{\lambda_a}$ are consistent regardless of the choice $a\in \{o,+,-\}$.
\begin{proposition}\label{prar}
    Consider $0<\alpha<2\pi$ and a monic elliptic operator $L_A$ with standard root $V$. For $\lambda \in \mathbb{C}\setminus \{0\}$, any solution $u_\lambda=r^\lambda v$ to $L_Au_\lambda=0$ is of the form (\ref{wunderdich}). Similarly, any solution to $\mathcal{L}_A(\lambda)v_\lambda=0$ is of the form (\ref{wunderdich2}).
\end{proposition}
\begin{proof}
By the preceding discussion and the final remark in Section \ref{eigenvalofop}, it suffices to show that $u_\lambda=0$ implies $0=(c_1,c_2)\in \mathbb{C}^\ell$. Let us assume the contrary and derive a contradiction. Then there are $c_1,c_2\in \mathbb{C}^\ell$, at least one $c_\bullet\neq 0$, such that
    \begin{align}\label{diffi}
        (\cos(\varphi) \Id_\ell+\sin(\varphi) V)^{\lambda_+}c_1+(\cos(\varphi) \Id_\ell +\sin(\varphi)\overline{V})^{\lambda_-} c_2=0\quad \text{for } \varphi\in [0,\alpha).
    \end{align}
In particular, for $\varphi=0$, this leads to $-c_2=c_1\neq 0$. Differentiating (\ref{diffi}) and evaluating at $\varphi=0_+$, we get the condition $(V-\overline{V})c_1=0$, which implies $\det(\Im V)=0$, contradicting Theorem \ref{algtheo}.
\end{proof}

\begin{remark}
i) The main reason for introducing monic reductions and standard roots was precisely to obtain explicit formulas such as (\ref{wunderdich}). This works because terms of the form $(x_1\Id_\ell + x_2 V)^\lambda$ can be handled via the functional calculus. By contrast, expressions like $(x_1 V_1+x_2 V_2)^\lambda$ with non-commuting $V_1,V_2\in \Mat_\ell(\mathbb{C})$ cannot be treated in the same way, and their spectral properties and derivatives are not accessible without further structure.\\ 
ii) Note that we excluded $\lambda=0$ because $u_0(x_1,x_2)=c_1+c_2$ for $c_1,c_2\in \mathbb{C}^\ell$ does not yield $2\ell$ linearly independent solutions. This is why the case $\lambda=0$ was treated separately in Lemma \ref{lambda=0}.\\
iii) The representation of solutions for the model problem without boundary conditions provided in (\ref{wunderdich}) is new. However, another representation involving complex contour integrals is known and can be found in \cite{CoDauCCSA1993} and \cite{CoDauFSAC2001}.\\
iv) The reason for avoiding $\alpha=2\pi$ in the discussion and Prop. \ref{prar} is to prevent ill-defined $\arg_a$-functions. Nonetheless, Prop. \ref{prar} can be canonically generalized to $\alpha=2\pi$ by considering a continuous continuation of $\arg_a$ from the case $\alpha<2\pi$.
\end{remark}

\section{Analysis of the model problem with boundary conditions}\label{analysis of the model problem}
After deriving explicit formulas for solutions of the model problem without boundary conditions, we now investigate solutions of the model problem with boundary conditions. It will be shown that existence of a solution $r^\lambda v$ for angle $\alpha$ is equivalent to $0\in \sigma(M_{\lambda,\alpha})$ for some matrix $M_{\lambda,\alpha}\in \Mat_\ell(\mathbb{C})$, which depends on the boundary condition.
\smallskip

For the remainder of this section, assume $0<\alpha< 2\pi$, and a monic elliptic tuple $A=(A_{11},A_{12},\Id_\ell)$ with standard root $V$. We consistently write $V=(S+i\Id_\ell)D=C+iD$ for $C,D,S\in \Mat_\ell(\mathbb{R})$. By Prop. \ref{prar}, we have that $u=r^\lambda v$ for $\lambda \in \mathbb{C}\setminus \{0\}$ is a solution to $ L_A(\partial_{x_1},\partial_{x_2}) u=0$ if there exist $c_1,c_2\in \mathbb{C}^\ell$ such that:
\begin{align}\label{solutionv}
    v(\varphi)=(\cos(\varphi) \Id_\ell+\sin(\varphi) V)^{\lambda_+}c_1+(\cos(\varphi) \Id_\ell +\sin(\varphi)\overline{V})^{\lambda_-} c_2.
\end{align}
\subsection{Dirichlet boundary conditions}
Assume that $r^\lambda v$, for $v$ in (\ref{solutionv}), satisfies Dirichlet boundary conditions on $\Gamma^\pm$. Consequently,
\begin{align}\label{newlab}
    0=v(0)=c_1+c_2,\quad 0=v(\alpha)=V_\alpha^{\lambda_+}c_1+\overline{V_\alpha}^{\lambda_-} c_2,
\end{align}
where we denote 
\begin{align}\label{VAA}
    V_\alpha=\cos(\alpha)\Id_\ell+ \sin(\alpha)V=\cos(\alpha)\Id_\ell+\sin(\alpha)SD +i\sin(\alpha)D.
\end{align}
Rewriting (\ref{newlab}), we have $(V_\alpha^{\lambda_+}-\overline{V_\alpha}^{\lambda_-})c_1=0 $. Thus, finding a solution $u=r^\lambda v\neq 0$ for the model problem with Dirichlet boundary conditions is equivalent to $0\in \sigma(V_\alpha^{\lambda_+}-\overline{V_\alpha}^{\lambda_-})$. Observe that:
\begin{align*}
     0\in \sigma(V_\alpha^{\lambda_+}-\overline{V_\alpha}^{\lambda_-})\iff 0\in\sigma (D^{1/2}V_\alpha^{\lambda_+} D^{-1/2}-D^{1/2}\overline{V_\alpha}^{\lambda_+} D^{-1/2}),
\end{align*}
using $D>0$ (Lemma \ref{inv}). We set
$Z_\alpha:=D^{1/2}V_\alpha D^{-1/2}$ and note that
\begin{align}\label{Za}
   Z_\alpha=\cos(\alpha)\Id_\ell+D^{1/2}SD^{1/2} \sin(\alpha)+iD\sin(\alpha)
\end{align}
By properties of the functional calculus, we have $D^{1/2}V_\alpha^{\lambda_+} D^{-1/2} =Z_\alpha^{\lambda_+}$ (Appendix \ref{spundfunc}). The discussion shows the following result.
\begin{proposition}\label{dirstate}
    Consider $0<\alpha<2\pi$ and a monic elliptic operator $L_A$ with standard root $V=(S+i\Id_\ell)D$. The model problem with Dirichlet boundary conditions has for $\lambda \in \mathbb{C}\setminus \{0\}$ a solution $r^\lambda v\neq 0$ if and only if
    \begin{align}\label{thedireq}
        0\in \sigma(M_{\lambda,\alpha})\quad \text{for }M_{\lambda,\alpha}=Z_\alpha^{\lambda_+}-\overline{Z_\alpha}^{\lambda_-}.
\end{align}
\end{proposition}

\subsection{Mixed boundary conditions}

First, we characterize the Neumann boundary condition $\mathcal{N}_A v(\varphi)=0$ for $\varphi\in \{0,\alpha\}$. For this, calculate (note that $n=(-\sin(\varphi),\cos(\varphi))=\frac{1}{r}(-x_2,x_1) $ is the normal vector):
\begin{align*}
    \sum_{i,j=1}^2 A_{ij} n_i \partial_{x_j}(x_1 \Id_\ell+x_2 V)^{\lambda_+}= \frac{\lambda}{r} \big(-x_2(A_{11}+A_{12}V)+x_1(A_{12}+ V) \big) (x_1\Id_\ell+x_2 V)^{(\lambda-1)_+}.
\end{align*}
Since $V$ solves $A_{11}+A_{12}V=- A_{12}V - V^2$ (Theorem \ref{algtheo}), we obtain
\begin{align*}
  \sum_{i,j=1}^2 A_{ij} n_i \partial_{x_j}(x_1 \Id_\ell+x_2 V)^{\lambda_+}=\frac{2\lambda}{r}   (A_{12}+V)(x_1\Id_\ell+x_2 V)^{\lambda_+}.
\end{align*}
Similarly, we derive $ \sum_{i,j=1}^2 A_{ij} n_i \partial_{x_j}(x_1 \Id_\ell+x_2 \overline{V})^{\lambda_-}=\frac{2\lambda}{r}   (A_{12}+\overline{V})(x_1\Id_\ell+x_2 \overline{V})^{\lambda_-}$.
\smallskip
Now, assume that $r^\lambda v\neq 0$, for $v$ in (\ref{solutionv}), satisfies Dirichlet boundary conditions on $\Gamma^+$ and Neumann boundary conditions on $\Gamma^-$. Thus, the coefficients $c_\bullet$ in (\ref{solutionv}) satisfy (recall that we assumed $\lambda \neq 0$):
\begin{align}\label{lupau}
     0=c_1+c_2,\quad
    0=(A_{12}+V)V_\alpha^{\lambda_+}c_1 +(A_{12}+\overline{V})\overline{V_\alpha}^{\lambda_-}c_2,
\end{align}
for $V_\alpha$ in (\ref{VAA}). This is equivalent to $0\in \sigma(\widetilde{M}_{\lambda,\alpha})$ for
\begin{align*}
    \widetilde{M}_{\lambda,\alpha}=(A_{12}+V)V_\alpha^{\lambda_+} -(A_{12}+\overline{V})\overline{V_\alpha}^{\lambda_-}.
\end{align*}
Write $ A_{12}=-\frac{1}{2}(V+V^*)$ and $V=(S+i\Id_\ell)D$ (see Theorem \ref{algtheo}) such that:
\begin{align}\label{pumpatt}
    A_{12}+V=\frac{1}{2}[S,D]+iD,\quad A_{12}+\overline{V}=\frac{1}{2}[S,D]-iD,
\end{align}
and we conclude:
\begin{align*}
    0\in \sigma (\widetilde{M}_{\lambda,\alpha})\iff
    0\in \sigma\left(\left(\frac{1}{2}[S,D]+iD\right)V_\alpha^{\lambda_+} -\left(\frac{1}{2}[S,D]-iD\right)\overline{V_\alpha}^{\lambda_-} \right).
\end{align*}

Now, again using $Z_\alpha=D^{1/2}V_\alpha D^{-1/2}$, this is equivalent to:
\begin{align*}
    0\in~&\sigma\left( D^{-1/2}\left(\frac{1}{2}[S,D]+iD\right) D^{-1/2}Z_\alpha^{\lambda_+}-D^{-1/2}\left(\frac{1}{2}[S,D]-iD\right)D^{-1/2} \overline{Z_\alpha}^{\lambda_-}\right)\\
    =~&\sigma\left( \frac{1}{2} D^{-1/2}[S,D]D^{-1/2}(Z_\alpha^{\lambda_+}-\overline{Z_\alpha}^{\lambda_-}) +i(Z_\alpha^{\lambda_+}+\overline{Z_\alpha}^{\lambda_-})\right).
\end{align*}
This leads to the following result.
\begin{proposition}\label{mixstate}
    Consider $0<\alpha<2\pi$ and a monic elliptic operator $L_A$ with standard root $V=(S+i\Id_\ell)D$. The model problem with mixed boundary conditions has for $\lambda \in \mathbb{C}\setminus \{0\}$ a solution $r^\lambda v\neq 0$ if and only if
    \begin{align*}
        0\in \sigma\left(M_{\lambda,\alpha}\right)\quad \text{for }M_{\lambda,\alpha}= \frac{1}{2} [D^{-1/2}SD^{-1/2},D](Z_\alpha^{\lambda_+}-\overline{Z_\alpha}^{\lambda_-}) +i(Z_\alpha^{\lambda_+}+\overline{Z_\alpha}^{\lambda_-}).
\end{align*}
\end{proposition}

\subsection{Neumann boundary conditions}\label{neumann boundary yea}
So far, we have only set ellipticity conditions for the operator $L_A$ (Def. \ref{the conda}) but no conditions for $B^\pm_A$. For the model problem to form an elliptic system in the sense of Agmon-Douglis-Nirenberg, it is required that $B^\pm_A$ satisfies the so called \textit{complementing boundary condition}, which is discussed in more detail in Appendix \ref{ellli}. A discussion of the complementing boundary condition in the case of Dirichlet boundary conditions for $\Gamma^-$ was not necessary, as it is automatically satisfied for strongly elliptic systems (Remark 3.2.7 in \cite{CDN2010CSAR}), however, it is necessary for Neumann boundary conditions.
\smallskip

To motivate the following definitions, note that the conditions on $c_1$, $c_2$ for Neumann boundary conditions on $\Gamma^\pm$ have the form:
\begin{align}\label{neumannni}
     0=(A_{12}+V)c_1+(A_{12}+\overline{V})c_2,\quad
    0=(A_{12}+V)V_\alpha^{\lambda_+}c_1 +(A_{12}+\overline{V})\overline{V_\alpha}^{\lambda_-}c_2,
\end{align}
which is deduced analogously to (\ref{lupau}). Solving for $c_1$ in the first equation requires $A_{12}+V$ to be invertible, which is equivalent to the \textit{complementing boundary condition} for $N_A^\pm$ (see Appendix \ref{ellli}).
\begin{definition}\label{neumann and contr neumann}
    Consider an elliptic tuple $A$, where $V=C+iD$ for $C,D\in \Mat_\ell(\mathbb{R})$ is the standard root. We say
    \begin{enumerate}
        \item[i)] $A$ is \textit{Neumann well-posed} if $2i\notin \sigma([D^{-1},C])$.
        \item[ii)] $A$ is \textit{contractive Neumann well-posed} if $\rho([D^{-1},C])< 2$.
    \end{enumerate}
\end{definition}
How do these definitions relate to the considerations discussed above? From (\ref{pumpatt}):
\begin{align}\label{byebye}
    A_{12}+V=\frac{1}{2}[S,D]+iD=\frac{1}{2}D(D^{-1}SD-S+2i\Id_\ell),
\end{align}
which is invertible if and only if $-2i$ is not an eigenvalue of:
\begin{align}\label{recalli}
    D^{-1}SD-S=D^{-1}C-CD^{-1}=[D^{-1},C].
\end{align}
Note that $\sigma([D^{-1},C])\subset i \mathbb{R}$ and $\overline{\sigma([D^{-1},C])}= \sigma([D^{-1},C])$. This follows from:
\begin{align}\label{badimbadum}
    \sigma([D^{-1},C])=\sigma(D^{1/2}[D^{-1},C]D^{-1/2})\stackrel{(\ref{recalli})}{=}\sigma([D^{-1/2}SD^{-1/2},D])\subset i \mathbb{R},
\end{align}
where the last inclusion holds since $D^{-1/2}SD^{-1/2}$ and $D$ are symmetric, their commutator is skew-symmetric, and skew symmetric matrices have imaginary, complex conjugated eigenvalues. Thus, $-2i\notin \sigma([D^{-1},C])$ if and only if $2i\notin \sigma([D^{-1},C])$. So, by (\ref{byebye}), we have shown that invertibility of $A_{12}+V$ is equivalent to Neumann well-posedness. Contractive Neumann well-posedness will become relevant for the main result Theorem \ref{main mixed}

\begin{remark}
    Assuming that $A$ is an elliptic tuple, the various ellipticity conditions discussed in this work are related as follows:
    \begin{align*}
        \text{A formal positive}&\implies \text{A contractive Neumann well-posed}
        \implies \text{A Neumann well-posed}\\
        &\iff \text{A fulfills complementing boundary conditions for $N_A$}.
    \end{align*}
    Here, formal positivity is a common ellipticity condition in linear elasticity. Its definition is given in Appendix \ref{ellli}, where one can also find proofs for the above implications
\end{remark}

We proceed with Neumann boundary conditions for the model problem. Assume that $r^\lambda v\neq 0$, for $v$ in (\ref{solutionv}), satisfies Neumann boundary conditions on $\Gamma^\pm$. \textit{Additionally, assume that the elliptic tuple $A$ is Neumann well-posed}. From (\ref{neumannni}), we derive $0\in \sigma(\widetilde{M}_{\lambda,\alpha})$ for
\begin{align*}
    \widetilde{M}_{\lambda,\alpha}=(A_{12}+V)V_\alpha^{\lambda_+}(A_{12}+V)^{-1}-(A_{12}+\overline{V})\overline{V_\alpha}^{\lambda_-}(A_{12}+\overline{V})^{-1}.
\end{align*}
By substituting (\ref{pumpatt}), we get the condition:
\begin{align*}
   0\in \sigma \left(\left(\frac{1}{2}[S,D]+iD\right)V_\alpha^{\lambda_+} \left(\frac{1}{2}[S,D]+iD\right)^{-1}- \left(\frac{1}{2}[S,D]-iD\right)\overline{V_\alpha}^{\lambda_-} \left(\frac{1}{2}[S,D]-iD\right)^{-1}\right).
\end{align*}
Rewriting $Z_\alpha=D^{1/2}V_\alpha D^{-1/2}$ and using similar arguments as before, we arrive at the following result.
\begin{proposition}\label{neumannstate}
     Consider $0<\alpha<2\pi$ and a Neumann well-posed monic elliptic operator $L_A$ with standard root $V=(S+i\Id_\ell)D$. The model problem with Neumann boundary conditions has for $\lambda \in \mathbb{C}\setminus \{0\}$ a solution $r^\lambda v\neq 0$ if and only if $0\in \sigma(M_{\lambda,\alpha})$ for
    \begin{align*}
        M_{\lambda,\alpha}=E~ Z_\alpha^{\lambda_+} E^{-1}- \overline{E} ~\overline{Z_\alpha}^{\lambda_-} \overline{E}^{-1}\quad \text{where }E=\frac{1}{2}D^{-1/2}[S,D]D^{-1/2}+i\Id_\ell.
\end{align*}
\end{proposition}

\section{Investigation of the spectrum of $M_{\lambda,\alpha}$}
\label{matrix eq associated}

In the previous section, we derived the condition $0\in\sigma(M_{\lambda,\alpha})$ corresponding to solutions of the model problem subject to Dirichlet, mixed, or Neumann boundary conditions. This section provides bounds on $|\Re\lambda|$ that guarantee $0\notin \sigma(M_{\lambda})$ for matrices $M_\lambda\in \Mat_\ell(\mathbb{C})$ encoding the structure of $M_{\lambda,\alpha}$ in the Dirichlet and mixed case. These bounds will be crucial in deriving the main results in Section \ref{regularity results of the model problem}. Our strategy is to use the numerical range to bound the eigenvalues of $M_{\lambda}$ away from zero. 
\subsection{Numerical range}
This subsection provides a small introduction and main properties of the numerical range. The reference is §1 in \cite{HoJo1991TIMA}, where the numerical range is referred to as "field of values". Additional results can be found in \cite{GuDu1997NRFV}. For $Z\in \Mat_\ell(\mathbb{C})$, the \textit{numerical range of $Z$} is defined as:
\begin{align*}
    W(Z)=\left\{\langle x,Zx\rangle: x\in \mathbb{C}^\ell, \|x\|=1\right\},
\end{align*} 
and the \textit{angular field of $Z$} as:
\begin{align*}
    W'(Z)=\left\{\langle x,Zx\rangle:x\in \mathbb{C}^\ell \setminus \{0\}\right\}.
\end{align*}
Using the scaling $x\mapsto r x$ for $r>0$, it is clear that $W'(Z)$ consists of rays connecting the origin $0\in \mathbb{C}$ to points in $W(Z)$, and $W(Z)\subset W'(Z)$. For $Z,Y\in \Mat_\ell(\mathbb{C})$ and $\alpha,\beta \in \mathbb{C}$, the following hold:
\begin{enumerate}[label=\textbf{N\arabic*}:]
    \item $W(Z)$ is convex and compact.
    \item $W(\alpha Z+\beta \Id_\ell)=\alpha W(Z)+\{\beta\}$.
    \item The spectrum is bounded by the numerical range: $\sigma(Z)\subset W(Z)$.
    \item $W(Z)=W(U^*ZU)$ for any unitary $U\in \Mat_\ell(\mathbb{C})$ and $W'(Z)=W'(C^*ZC)$ for any $C\in \Mat_\ell(\mathbb{C})$.
    \item $W(Z+Y)\subset W(Z)+W(Y)$.
    \item $W(Z)$ is a line segment $[\alpha,\beta]\subset \mathbb{R}$ if and only if $Z$ is Hermitian. Then $\min \{\sigma(Z)\}=\alpha$ and $\max \{\sigma(Z)\}=\beta$.
    \item $W(Z)\subset \operatorname{RHP}$ if and only if $Z+Z^*>0$.
    \item $W(Z^*)=\overline{W(Z)}$.
\end{enumerate}

The following definition encodes the structure of $Z_\alpha$ in (\ref{Za}) for $0<\alpha<\pi$.

\begin{definition}\label{def+i}
    The set of symmetric matrices $Z\in \Mat_\ell(\mathbb{C})$ satisfying $\Im(Z)>0$ is denoted by $\Mat_\ell(\mathbb{C})_{+i}$. 
\end{definition}

\begin{lemma}\label{symm}
Consider $Z\in \Mat_\ell(\mathbb{C})$ with $Z=Z^T$. Then $Z\in \Mat_\ell(\mathbb{C})_{+i}$ if and only if $W(Z)\subset \operatorname{UHP}$. In this case, $\sigma(Z)\subset \operatorname{UHP}$.
\end{lemma}
\begin{proof}
    The first statement follows from \textbf{N7} applied to $-i Z$, using $\Re Z=(\Re Z)^T$ and $\Im Z=(\Im Z)^T>0$. The second from \textbf{N3}.
\end{proof}

Two key results are presented which guide the subsequent proofs. The proofs of these statements are given in Appendix \ref{numerrange and accreive}.

\begin{lemma}\label{numrang}
     Consider $Z\in \Mat_\ell(\mathbb{C})_{+ i}$ and $\lambda\in [-1,1]\setminus \{0\}$. Then $W'(Z^\lambda)\subset  \operatorname{sgn}(\lambda) \operatorname{UHP}$.
\end{lemma}

\begin{lemma}\label{chaotisch}
    Consider $Z\in \Mat_\ell(\mathbb{C})_{+i}$ and $\lambda\in \mathbb{R}\setminus \{0\}$. Then the following holds:
    \begin{align*}
        &i) \quad \operatorname{sgn}(\lambda)\left(\Id_\ell-(Z^{i\lambda })^*Z^{i\lambda}\right)>0,
         &ii) \quad \min\{\beta: \beta \in \sigma((Z^{i\lambda })^*Z^{i\lambda})\}\xrightarrow[]{\lambda\to -\infty} \infty.
    \end{align*}
\end{lemma}
Note that, in this section, complex exponentiation, as well as $\arg$- and $\log$-functions, will always refer to the principal branch.
\subsection{Dirichlet boundary conditions}
\begin{theorem}\label{dir}
    Consider $Z\in \Mat_\ell(\mathbb{C})_{+i}$, $\lambda\in \mathbb{C}\setminus \{0\}$ with $|\Re \lambda|\leq 1 $, and define $M_\lambda=Z^{\lambda}-\overline{Z}^{\lambda}$. Then $0\notin \sigma(M_\lambda)$.
\end{theorem}
\begin{proof}
Decompose $\lambda$ into real and imaginary part $\lambda=\lambda_1+i\lambda_2$ and write:
\begin{align} 
        0\in &~\sigma(Z^\lambda-\overline{Z}^\lambda )=\sigma(Z^{i\lambda_2}Z^{\lambda_1}-\overline{Z}^{\lambda_1} \overline{Z}^{i\lambda_2} )\label{nase}\\ \nonumber
        \iff 0\in&~ \sigma(Z^{i\lambda_2}Z^{\lambda_1}\overline{Z}^{-i\lambda_2} -\overline{Z}^{\lambda_1}  )=\sigma(Z^{i\lambda_2}Z^{\lambda_1}(Z^{i\lambda_2})^* -\overline{Z}^{\lambda_1}  ).
\end{align}
Here, we used $Z^{i\lambda_2}Z^{\lambda_1}=Z^{\lambda_1}Z^{i\lambda_2}$, $\det \big(\overline{Z}^{i\lambda_2}\big)\neq 0$, and $(\overline{Z}^{i\lambda_2})^{-1}=\overline{Z}^{-i\lambda_2}=(Z^{i\lambda_2})^*$ (see Appendix \ref{spundfunc}).
\medskip

\Needspace{3\baselineskip}
\textbf{Case 1: }$\Re \lambda \neq 0$.

Let us assume $\lambda_1\in [-1,1]\setminus \{0\}$. We aim to show that the spectrum of $Z^{i\lambda_2}Z^{\lambda_1}(Z^{i\lambda_2})^* -\overline{Z}^{\lambda_1}  $ is bounded away from $0$, implying that $0\in \sigma(Z^\lambda-\overline{Z}^\lambda)$ is not possible. We derive:
\begin{align}\nonumber
    &\sigma( Z^{i\lambda_2} Z^{\lambda_1}(Z^{i\lambda_2})^* -\overline{Z}^{\lambda_1} )\subset W( Z^{i\lambda_2} Z^{\lambda_1}(Z^{i\lambda_2})^* -\overline{Z}^{\lambda_1} )\subset W( Z^{i\lambda_2} Z^{\lambda_1}(Z^{i\lambda_2})^*) -W(\overline{Z}^{\lambda_1} )\\ \label{aswehave}
    \subset ~&W'(  Z^{\lambda_1}) -\overline{W(Z^{\lambda_1} )}\subset \operatorname{sgn}(\lambda_1) \operatorname{UHP} -\operatorname{sgn}(\lambda_1) \overline{\operatorname{UHP}}\subset \operatorname{sgn}(\lambda_1) \operatorname{UHP},
\end{align}
by properties \textbf{N3}, \textbf{N5}, \textbf{N4}, \textbf{N8} of the numerical range, $\overline{Z}^{\lambda_1}=(Z^{\lambda_1})^*$, Lemma \ref{numrang}, and $\overline{\operatorname{UHP}}=-\operatorname{UHP}$. Note that $0\notin \operatorname{sgn}(\lambda_1)\operatorname{UHP}$.
\medskip

\textbf{Case 2: }$\Re \lambda = 0$.

Next, assume $\lambda_1=0$ and $\lambda_2\in \mathbb{R}\setminus \{0\}$. By (\ref{nase}), the argument boils down to show that $0\notin \sigma(Z^{i\lambda_2}(Z^{i\lambda_2})^* -\Id_\ell  )$. This follows from Lemma \ref{chaotisch}.
\end{proof}
\subsection{Mixed boundary conditions}
\begin{theorem}\label{mix}
Consider $Z\in \Mat_\ell(\mathbb{C})_{+i}$, $A,B\in \Mat_\ell(\mathbb{R})$ symmetric, $\lambda\in \mathbb{C}$ with $|\Re \lambda|\leq  \frac{1}{2}$, and define $M_\lambda= [A,B](Z^\lambda-\overline{Z}^\lambda)+i(Z^\lambda+\overline{Z}^\lambda)$. Then:
\begin{enumerate}
    \item[i)] $0\notin \sigma(M_\lambda)$ for $|\Re \lambda| \in (0,\frac{1}{2}]$.
    \item[ii)] $0\in \sigma(M_\lambda)$ for $\Re \lambda=0$ if and only if $ \rho([A,B])> 1$.
\end{enumerate}
\end{theorem}
\noindent
The condition $\rho([A,B])>1$ will later be linked to the failure of contractive Neumann well-posedness. For the proof, note that $[A,B]\in \Mat_\ell(\mathbb{R})$ is skew-symmetric due to symmetry of $A,B$. Thus, $\sigma([A,B])\subset i \mathbb{R}$, and all eigenvalues come in conjugate pairs. Similarly, $i[A,B]$ is Hermitian, so $\sigma(i[A,B])\subset \mathbb{R}$.

\begin{proof}
    Consider $Z$,$A$,$B$, and $\lambda\in \mathbb{C}$ as given in the statement. Decompose $\lambda$ into real and imaginary part $\lambda=\lambda_1+i\lambda_2$ and note, using the same arguments as in the derivation for (\ref{nase}):
\begin{align}\nonumber
     0\in&~\sigma \left( \big([A,B](Z^\lambda-\overline{Z}^\lambda)+i(Z^\lambda+\overline{Z}^\lambda)\right)\\ \label{whatyoulookingat}
     \iff 0\in&~\sigma\left( [A,B](Z^{i\lambda_2} Z^{\lambda_1}(Z^{i\lambda_2})^* -\overline{Z}^{\lambda_1}  )+i(Z^{i\lambda_2} Z^{\lambda_1}(Z^{i\lambda_2})^* +\overline{Z}^{\lambda_1}  )\right).
\end{align}

\textbf{Case 1: }$\lambda_1\neq 0$.

In this case, $R_\lambda:=Z^{i\lambda_2} Z^{\lambda_1}(Z^{i\lambda_2})^* -\overline{Z}^{\lambda_1} $ is invertible as shown in (\ref{aswehave}) since $0\notin \sigma(R_\lambda)$. Thus, due to (\ref{whatyoulookingat}) and invertibility of $R_\lambda^*$, $0\in \sigma(M_\lambda)$ is equivalent to:
\begin{align*}
    0\in &~\sigma \big( R_\lambda^* [A,B]R_\lambda+iR_\lambda^* \big(R_\lambda+2\overline{Z}^{\lambda_1}\big) \big).
\end{align*}
For Case 1, it suffices, by \textbf{N3}, to prove the following claim.
\medskip

\textbf{Claim 1: }$W \big( R_\lambda^* [A,B]R_\lambda+iR_\lambda^* \big(R_\lambda+2\overline{Z}^{\lambda_1}\big) \big) \subset \operatorname{sgn}(\lambda_1)\operatorname{RHP}$ for $\lambda_1\in [-\frac{1}{2},\frac{1}{2}]\setminus \{0\}$.

Note that:
\begin{align*}
    W(R_\lambda^*[A,B]R_\lambda)&\subset W'(R_\lambda^*[A,B]R_\lambda)\subset W'([A,B])\subset -i\cdot W'(i[A,B])\subset i \mathbb{R},\\
    W \big( iR_\lambda^*R_\lambda \big)&\subset iW(R_\lambda^*R_\lambda)\subset i\mathbb{R}.
\end{align*}
In the first line, we used \textbf{N4}, \textbf{N2}, \textbf{N6}, and $i[A,B]$ being Hermitian, and in the second line, we used \textbf{N2} and \textbf{N6}. By additivity of the numerical range \textbf{N5}, it suffices to show $W(iR_\lambda^*\overline{Z}^{\lambda_1})\subset \operatorname{sgn}(\lambda_1)\operatorname{RHP}$ for the claim. We have
\begin{align*}
    R_\lambda^*\overline{Z}^{\lambda_1}=Z^{i\lambda_2} \overline{Z}^{\lambda_1}(Z^{i\lambda_2})^*\overline{Z}^{\lambda_1} -Z^{\lambda_1}\overline{Z}^{\lambda_1}=Z^{i\lambda_2} \overline{Z}^{2\lambda_1}(Z^{i\lambda_2})^* -Z^{\lambda_1}(Z^{\lambda_1})^*,
\end{align*}
where we used $(Z^{i\lambda_2})^*\overline{Z}^{\lambda_1}=\overline{Z}^{\lambda_1}(Z^{i\lambda_2})^*$ and $\overline{Z}^{\lambda_1}=\overline{Z^{\lambda_1}}=(Z^{\lambda_1})^*$. We conclude
\begin{align*}
W(iR_\lambda^*\overline{Z}^{\lambda_1})\subset i W(Z^{i\lambda_2} \overline{Z}^{2\lambda_1}(Z^{i\lambda_2})^*)-iW(Z^{\lambda_1}(Z^{\lambda_1})^*)\subset i W'(\overline{Z}^{2\lambda_1})+i\mathbb{R}\subset \operatorname{sgn}(\lambda_1)\operatorname{RHP},
\end{align*}
by \textbf{N2}, \textbf{N5}, \textbf{N4}, \textbf{N6}, \textbf{N8}, Lemma \ref{numrang}, and $i\overline{\operatorname{UHP}}=\operatorname{RHP}$. This shows Claim 1 and closes Case 1.
\medskip

\textbf{Case 2: }$\lambda_1=0$.

In this case, (\ref{whatyoulookingat}) can be rewritten as:
\begin{align}\label{sofro}
    0\in \sigma\left( [A,B](Z^{it}(Z^{it})^* -\Id_\ell  )+i(Z^{it}(Z^{it})^* +\Id_\ell  )\right),
\end{align}
where we write $\lambda_2=t\in \mathbb{R}$ in the following. For $t=0$, this is equivalent to $0\in \sigma\left(2i\Id_\ell\right)$, which is not possible. So assume $t\neq 0$ such that $-\operatorname{sgn}(t) (Z^{it}(Z^{it})^* -\Id_\ell  )>0$ due to Lemma \ref{chaotisch}. Define for $t\in \mathbb{R}\setminus \{0\}$:
\begin{align*}
    K_t:=\left(Z^{it} (Z^{it})^*-\Id_\ell\right)\left(Z^{it} (Z^{it})^*+\Id_\ell\right)^{-1}.
\end{align*}
Alternatively, $K_t$ can be given by the functional calculus $K_t=f(Z^{it} (Z^{it})^*)$ for 
\begin{align}\label{effie}
    f:\mathbb{C}\setminus \{-1\}\to \mathbb{C},\quad z\mapsto \frac{z-1}{z+1}.
\end{align}
Note that $K_t$ is Hermitian, and by the spectral mapping theorem (Appendix \ref{spundfunc}), along with Lemma \ref{chaotisch}, one has $-\operatorname{sgn}(t)K_t>0$. Now (\ref{sofro}) is equivalent to:
\begin{align*}
    (\ref{sofro})&\iff 0\in \sigma([A,B]K_t+i \Id_\ell)\iff -i\in \sigma( [A,B] K_t)\\
    &\iff i \operatorname{sgn}(t)\in \sigma(-[A,B] \operatorname{sgn}(t)K_t) \iff i \operatorname{sgn}(t)\in  \sigma(|K_t|^{1/2}[A,B]|K_t|^{1/2}),
\end{align*}
where we used the abbreviation $|K_t|=-\operatorname{sgn}(t)K_t$. Before we continue, let us show the following claim.
\medskip

\textbf{Claim 2: }$\rho(|K_t|^{1/2}[A,B]|K_t|^{1/2})<\rho([A,B])$ for any $t\in \mathbb{R}$.

By the spectral mapping theorem, (\ref{effie}), and $\sigma(Z^{it} (Z^{it})^*)\subset \mathbb{R}_{>0}$, we deduce that $\rho(K_t)=\| K_t \|<1$, since $K_t$ is Hermitian. Consequently, we have $\|[A,B] |K_t| \| < \|[A,B]\|=\rho([A,B])$. The claim is deduced by:
\begin{align*}
    \rho(|K_t|^{1/2}[A,B]|K_t|^{1/2})=\rho([A,B]|K_t|)\leq \| [A,B]|K_t|\|<\rho([A,B]).
\end{align*}

It remains to show $0\in \sigma(M_{it})$ for some $t \in \mathbb{R}$ if and only if $\rho[A,B] )>1$. So far, we have shown $0\in \sigma(M_{it})$ if and only if
\begin{align}\label{ograe}
     t\neq0 \text{ and }i\operatorname{sgn}(t)\in  \sigma(|K_t|^{1/2}[A,B]|K_t|^{1/2}).
\end{align}
    
\textbf{Case 2a: }$\rho([A,B] )\leq 1$.

If $\rho([A,B])\leq 1$, then Claim 2 implies $\rho(|K_t|^{1/2}[A,B]|K_t|^{1/2})<1$ for any $t\in \mathbb{R}\setminus \{0\}$. In this case, (\ref{ograe}) cannot hold, and thus, $0\notin \sigma(M_{it})$.
\medskip

\textbf{Case 2b: }$\rho([A,B] )>1$.

In this case, there exists some $\beta\in \sigma (i[A,B])$ such that $\beta>1$. To complete the proof, it suffices to consider from now on $t<0$ such that $|K_t|=K_t$. Due to Lemma \ref{chaotisch} and the definition of $K_t$, we have $K_t\xrightarrow[]{t\to 0_{-}} 0$ and $K_t\xrightarrow[]{t\to -\infty} \Id_\ell$ such that:
\begin{align*}
    i K^{1/2}_t[A,B]K^{1/2}_t \xrightarrow[]{t\to 0_{-}} 0,\quad  i K^{1/2}_t[A,B]K^{1/2}_t \xrightarrow[]{t\to -\infty} i[A,B].
\end{align*}
Note that $i K^{1/2}_t[A,B]K^{1/2}_t$ is Hermitian and admits only real eigenvalues. By continuity of eigenvalues, and since $\beta>1$, there exist some $t<0 $ such that $1\in \sigma(i K_t^{1/2}[A,B]K_t^{1/2})$. Consequently, $-i\in \sigma(K_t^{1/2}[A,B]K_t^{1/2})$, which shows that (\ref{ograe}), and thus $0\in \sigma(M_{it})$, holds.
\end{proof}
The next result is needed for mixed boundary conditions in the case where $\alpha\in \{\pi,2\pi\}$.

\begin{theorem}\label{mix2}
Consider $Z\in \Mat_\ell(\mathbb{C})_{+i}$, $A,B\in \Mat_\ell(\mathbb{R})$ symmetric, $\lambda\in \ \mathbb{C}$, and define 
\begin{align*}
    M_\lambda= [A,B](e^{i\lambda \pi}-e^{-i\lambda \pi})+i(e^{i\lambda \pi}+e^{-i\lambda \pi})\Id_\ell.
\end{align*}
Then:
\begin{enumerate}
    \item[i)] $0\in \sigma(M_\lambda)$ implies $\Re \lambda \in \frac{1}{2}\mathbb{Z}$.
    \item[ii)] If $\rho([A,B])\leq 1$, then $0\in \sigma(M_\lambda)$ implies $\Re \lambda =\frac{1}{2}+\mathbb{Z}$.
    \item[iii)] If $\rho([A,B])<1$, then for any $k\in \mathbb{Z}$ there exist $\ell$ solutions $\lambda \in \mathbb{C}$ with $\Re \lambda =\frac{1}{2}+k$ satisfying $0\in \sigma(M_\lambda)$ (counted with algebraic multiplicity).
\end{enumerate}
\end{theorem}

\begin{proof}
i) Note that $0\in \sigma(M_\lambda)$ is equivalent to
\begin{align}\label{sincos}
     0\in \sigma \left( [A,B]\sin(\lambda \pi)+\cos(\lambda \pi)\Id_\ell\right).
\end{align}
Observe that $\lambda \in \mathbb{Z}$ cannot be a solution to (\ref{sincos}) since then $\sin(\lambda \pi)=0$ and $\cos(\lambda \pi)\neq 0$. Thus, we can assume $\sin(\lambda \pi)\neq 0$. Dividing by $\sin(\lambda \pi)$, the conditions becomes:
\begin{align}\label{suchalamb}
    0\in \sigma(M_\lambda) \iff -\frac{1}{\tan(\lambda \pi)}\in \sigma([A,B]).
\end{align}
Now assume $\lambda \in \mathbb{C}\setminus \mathbb{Z}$ solves the RHS in (\ref{suchalamb}). Since $\sigma([A,B])\in i\mathbb{R}$, it follows that $\tan(\lambda \pi)\in i\mathbb{R}$. Using the following representation of complex tangent:
 \begin{align*}
        \tan(x+iy)=\frac{\sin(2x)+i\sinh(2y)}{\cosh(2y)+\cos(2x)}\quad \text{for }x,y\in \mathbb{R},
    \end{align*}
we deduce that $0\in \sigma(M_\lambda)$ implies $0=\sin(2\pi \Re \lambda)$, equivalent to $\Re \lambda \in \frac{1}{2}\mathbb{Z}$.

ii) Assume $\rho([A,B])\leq 1$. Write $\lambda=\frac{k}{2}+it$ for $k\in \mathbb{Z}$ and $t\in \mathbb{R}$. Using the tangent representation, we have:
\begin{align}\label{lohlaf}
    \frac{1}{\tan(\lambda \pi)}=-i\frac{\cosh(2t\pi )+(-1)^{k}}{\sinh(2t\pi )}.
\end{align}

If $k\in 2 \mathbb{Z}$, then $\left| \tan(\lambda \pi)^{-1}\right|>1$. In this case, $\lambda$ cannot solve (\ref{sincos}) due to (\ref{suchalamb}) and $\rho([A,B])\leq 1$.

iii) Assume $\rho([A,B])<1$, so $\sigma([A,B])\subset i(-1,1)$. If $k\in 1+2\mathbb{Z}$, then the RHS of (\ref{lohlaf}) can be continuously extended to $t=0$ (with the value $0$) such that it defines a surjective function $f:\mathbb{R}\to i(-1,1)$. By (\ref{suchalamb}) and $\sigma([A,B])\subset i(-1,1)$, $0\in \sigma(M_{\frac{k}{2}+it})$ has $\ell$ solutions (counted with multiplicity) as $t$ in (\ref{lohlaf}) varies over $(-\infty,\infty)$.
\end{proof}

\section{Regularity results for the model problem}
\label{regularity results of the model problem}
The main results of this work, bounds on $|\Re \lambda|$ for Dirichlet and mixed boundary conditions, are given. Neumann boundary conditions are briefly discussed in Section \ref{neumanni}.
\subsection{Dirichlet boundary conditions}
The result in this case is not new, compare with §8.6 and §11.3 in \cite{VMR2001SPCS}. However, the proof is new and utilizes the methods established in this work.
\begin{theorem}\label{main dir}
    Consider an elliptic tuple $A=(A_{11},A_{12},A_{22})$. Define
    \begin{align*}
        \Lambda_\alpha:=\{\lambda \in \mathbb{C}: \exists\, r^\lambda v\neq 0 \text{ solving (\ref{soe}) with angle $\alpha$ and Dirichlet b.c.}\}
    \end{align*}
    Then:
   \begin{enumerate}[label=\roman*)]
    \item $\Lambda_\alpha \subset \{\lambda \in \mathbb{C}:|\Re \lambda| > 1\}~$ for $~0<\alpha<\pi$.
    \item $\Lambda_\pi =\mathbb{Z}\setminus \{0\}$.
    \item $\Lambda_\alpha \subset \{\lambda \in \mathbb{C}:|\Re \lambda| > \frac{1}{2}\}~$ for $~\pi<\alpha<2\pi$.
    \item $\Lambda_{2\pi} =\frac{1}{2}\mathbb{Z}\setminus \{0\}$.
\end{enumerate}
\end{theorem}
\begin{proof}
Due to Lemma \ref{lambda=0}, we can assume $\lambda\neq 0$ in the following. Furthermore, by Lemma \ref{reducioo}, it suffices to prove the result for monic elliptic tuples. For such tuples, Theorem \ref{algtheo} guarantees the existence of a standard root $V=(S+i \Id_\ell)D\in \Mat_\ell(\mathbb{C})$, where $S,D\in \Mat_\ell(\mathbb{R})$ are symmetric, and $D>0$. Due to Prop. \ref{dirstate}, the corresponding model problem with Dirichlet boundary conditions admits a solution $r^\lambda v\neq 0$ for $\lambda \in \mathbb{C}\setminus \{0\}$ and $0<\alpha<2\pi$ if and only if
\begin{align}\label{conti}
    0\in \sigma( Z_\alpha^{\lambda_+}-\overline{Z_\alpha}^{\lambda_-}).
\end{align}
Here, $ Z_\alpha=\cos(\alpha)+D^{1/2}SD^{1/2} \sin(\alpha)+iD\sin(\alpha)$ is a symmetric matrix. We now analyze the different cases for $\alpha$:
\smallskip

i) $0<\alpha<\pi$.

In this case, $Z_\alpha\in \Mat_\ell(\mathbb{C})_{+i}$ due to $D>0$ and $\sin(\alpha)>0$. By Lemma \ref{symm}, this implies $\sigma(Z_\alpha)\subset \operatorname{UHP}$ and $\sigma(\overline{Z_\alpha})\subset -\operatorname{UHP}$. Consequently, we can replace $\lambda_\pm$ with the principal branch $\lambda_o=\lambda$, and (\ref{conti}) simplifies to $0\in \sigma( Z_\alpha^{\lambda}-\overline{Z_\alpha}^{\lambda})$. The results follows by Theorem \ref{dir} and  $\lambda \neq 0$.
\smallskip

ii) $\alpha=\pi$.

In this case, $Z_\pi=-\Id_\ell$, so (\ref{conti}) reads
\begin{align*}
    0\in \sigma(((-1)^{\lambda_+}-(-1)^{\lambda_-} )\Id_\ell)=\sigma((e^{i\pi\lambda}-e^{-i\pi\lambda})\Id_\ell )\iff \sin(\lambda \pi)=0.
\end{align*}
The last equation holds if and only if $\lambda \in \mathbb{Z}$. This shows $\Lambda_\pi=\mathbb{Z}\setminus \{0\}$ due to $\lambda \neq 0$.
\smallskip

iii) $\pi<\alpha<2\pi$.

Before we can reduce $\lambda_\pm$ to $\lambda$ as in i), we need a trick. Since $D>0$ and $\sin(\alpha)<0$, Lemma \ref{symm} implies $W(\overline{Z_\alpha})\subset \operatorname{UHP}$. Using Lemma \ref{numrang}, \textbf{N8}, and $\overline{\operatorname{UHP}}=-\operatorname{UHP}$, this shows $W(Z_\alpha^{1/2})\subset -\operatorname{UHP}$. Define $Y_\alpha=-Z_\alpha^{1/2}$, and observe $W(Y_\alpha)\subset \operatorname{UHP}$ and $Y_\alpha^2=Z_\alpha$. Since $\sigma(Y_\alpha)\subset \operatorname{UHP}$ by \textbf{N3}, it follows that $Y_\alpha^{2\lambda}=Y_\alpha^{(2\lambda)_+}=Z^{\lambda_+}_\alpha$. Similarly, we obtain $\overline{Y_\alpha}^{2\lambda}=\overline{Z_\alpha}^{\lambda_{-}}$. Thus, (\ref{conti}) reads $ 0\in \sigma( Y_\alpha^{2\lambda}-\overline{Y_\alpha}^{2\lambda})$. Note that $Y_\alpha\in \Mat_\ell(\mathbb{C})_{+i}$ due to $W(Y_\alpha)\subset \operatorname{UHP}$ and Lemma \ref{symm}. The result follows by Theorem \ref{dir} and $\lambda \neq 0$.
\smallskip

iv) $\alpha=2\pi$.

Here, $Z_{2\pi}=\Id_\ell$, and we cannot use the original definition of $\lambda_\pm$. However, as pointed out in the last remark of \Cref{explisolf}, the result can be obtained as a limit case of $\alpha<2\pi$. For this, define $Z_\varepsilon:=(1-i\varepsilon)\Id_\ell$ for $\varepsilon>0$. Note $Z_\varepsilon\xrightarrow[]{\varepsilon\to 0_+}Z_{2\pi}$, as well as
\begin{align*}
    (Z_\varepsilon)^{\lambda_+}\to e^{2\pi i \lambda}Id_\ell,\quad (\overline{Z}_\varepsilon)^{\lambda_-}\to e^{-2\pi i\lambda}Id_\ell.
\end{align*}
Thus, taking the limit, (\ref{conti}) becomes $0\in \sigma((e^{2i\pi\lambda}-e^{-2i\pi\lambda})\Id_\ell )$ equivalent to $\sin(2\lambda \pi)=0$. This shows $\Lambda_{2\pi}=\frac{1}{2}\mathbb{Z}\setminus \{0\}$.
\end{proof}

\subsection{Mixed boundary conditions}

\begin{theorem}\label{main mixed}
Consider an elliptic tuple $A=(A_{11},A_{12},A_{22})$. Define:
    \begin{align*}
        \Lambda_\alpha:=\{\lambda \in \mathbb{C}: \exists\, r^\lambda v\neq 0\text{ solving (\ref{soe}) with angle $\alpha$ and mixed b.c.}\}
    \end{align*}
\begin{enumerate}
    \item If $A$ is contractive Neumann well-posed, then:
   \begin{enumerate}[label=\roman*)]
    \item $\Lambda_\alpha ~\subset \{\lambda \in \mathbb{C}:|\Re \lambda| > \frac{1}{2}\}~$ for $~0<\alpha<\pi$.
    \item $\Lambda_\pi ~\subset \{\lambda\in \mathbb{C}:\Re \lambda \in \frac{1}{2}+\mathbb{Z}\}$.
    \item $\Lambda_\alpha ~\subset \{\lambda \in \mathbb{C}:|\Re \lambda| > \frac{1}{4}\}~$ for $~\pi<\alpha<2\pi$.
    \item $\Lambda_{2\pi} \subset\{\lambda\in \mathbb{C}:\Re \lambda \in \frac{1}{4}+\frac{1}{2}\mathbb{Z}\}$.
\end{enumerate}
\item If the assumption on contractive Neumann well-posedness is dropped, solutions not satisfying the above conditions must be of the form:
\begin{enumerate}
    \item $\Re \lambda=0~$ for i) and iii).
    \item $\Re \lambda \in \frac{1}{2}\mathbb{Z}~$ for ii).
    \item $\Re \lambda \in \frac{1}{4}\mathbb{Z}~$ for iv).
\end{enumerate}
\end{enumerate}

\end{theorem}
This result is new and was the central motivation for developing the framework introduced in this work. Many arguments in the proof are similar to those for Dirichlet boundary conditions, so we omit details where appropriate.
\begin{proof}
By Lemma \ref{lambda=0}, we can assume $\lambda\neq 0$ in the following. Furthermore, by Lemma \ref{reducioo} (and the definition of contractive Neumann well-posedness), it suffices to prove the result for monic elliptic tuples (which are contractive Neumann well-posed). For such tuples, let $V=(S+i \Id_\ell)D$, where $S,D\in \Mat_\ell(\mathbb{R})$ are symmetric and $D>0$, denote its standard root. Recall that contractive Neumann well-posedness means:
\begin{align}\label{numanwell}
    \rho([(\Im V)^{-1},\Re V])<2\stackrel{(\ref{badimbadum})}{\iff} \rho\left(\frac{1}{2} [D^{-1/2}SD^{-1/2},D]\right)<1.
\end{align}

Due to Prop. \ref{mixstate}, the model problem with mixed boundary conditions admits a solution $r^\lambda v\neq 0$ for $\lambda \in \mathbb{C}\setminus \{0\}$ and $0<\alpha<2\pi$ if and only if for $Z_\alpha$ given in (\ref{Za}):
    \begin{align}\label{mixi}
        0\in \sigma\left( \frac{1}{2} [D^{-1/2}SD^{-1/2},D](Z_\alpha^{\lambda_+}-\overline{Z_\alpha}^{\lambda_-}) +i(Z_\alpha^{\lambda_+}+\overline{Z_\alpha}^{\lambda_-})\right)
\end{align}

i) $0<\alpha<\pi$.

As in the case of Dirichlet boundary conditions, one can replace $\lambda_\pm=\lambda$ in (\ref{mixi}). The result for i) and (a) follows from Theorem \ref{mix} and (\ref{numanwell}).
\smallskip

ii) $\alpha=\pi$.

Similar to the case of Dirichlet boundary conditions, one can write $Z_\alpha^{\lambda_+}=e^{i\lambda \pi}$ and $\overline{Z_\alpha}^{\lambda_-}=e^{-i\lambda \pi}$ in (\ref{mixi}). The result for ii) and (b) follows from Theorem \ref{mix2} and (\ref{numanwell}).
\smallskip

iii) $\pi<\alpha<2\pi$.

By the same argument as in the Dirichlet case, we can rewrite (\ref{mixi}) as 
\begin{align*}
     0=\det\left( \frac{1}{2} [D^{-1/2}SD^{-1/2},D](Y_\alpha^{2\lambda}-\overline{Y_\alpha}^{2\lambda}) +i(Y_\alpha^{2\lambda}+\overline{Y_\alpha}^{2\lambda})\right),
\end{align*}
for some $Y_\alpha\in \Mat_\ell(\mathbb{C})_{+i}$. The result for iii) and (a) follows from Theorem \ref{mix} and (\ref{numanwell}).
\smallskip

iv)  $\alpha=2\pi$.

By the same limiting argument as in the Dirichlet case, we can write $Z_\alpha^{\lambda_+}=e^{2i\lambda \pi}$ and $\overline{Z_\alpha}^{\lambda_-}=e^{-2i\lambda \pi}$ in (\ref{mixi}). The result for iv) and (c) follows from Theorem \ref{mix2} and (\ref{numanwell}).
\end{proof}

Let us emphasize again that Theorem \ref{main mixed} is highly relevant in regularity theory since it provides lower bounds on $|\Re \lambda|$ for the model problem of an angle $0<\alpha\leq 2\pi$ for a broad class of elliptic systems. For further details, we refer to Theorem 8.1.7 and §8.3.1 in \cite{MaRO2010EEPD}.

\begin{remark}
i) For 1.) in Theorem \ref{main mixed}, one could in principle also include elliptic tuples where the standard root $V$ of the monic reduction satisfies $ \rho([(\Im V)^{-1},\Re V])=2$. This follows because the main ingredients, Theorem \ref{mix} and Theorem \ref{mix2}, cover this case. Then, however, the system is not Neumann well-posed.\\
ii) Disregarding the case mentioned above, contractive Neumann well-posedness precisely distinguishes scenarios in which purely imaginary solutions can occur - the key motivation for introducing its definition. Verifying contractive Neumann well-posedness can be challenging in practice, however, Appendix \ref{ellli} relates it to a stronger ellipticity condition that is commonly used in linear elasticity.
\end{remark}

\section{Summary and comments}\label{summary}
In this work, we analyzed the model problem for an elliptic system subject to Dirichlet, mixed, and Neumann boundary conditions within a new framework. For all three boundary conditions, we derived a condition of the form $0\in \sigma (M_{\lambda,\alpha})$, which characterizes the pairs $(\alpha,\lambda)$ such that the model problem with angle $\alpha$ admits a solution of the form $r^\lambda v$. For Dirichlet and mixed boundary conditions, we established lower bounds on $|\Re \lambda|$ for nontrivial solutions. For the former, these results align with those found in the literature. For the latter, our findings represent an important contribution for the regularity of linear second-order elliptic systems with mixed boundary conditions.

\subsection{Bounds on $|\Re \lambda|$ for Neumann boundary conditions}\label{neumanni}
\begin{figure}[h!]
    \centering
    \begin{subfigure}[b]{0.48\linewidth}
        \centering
        \includegraphics[width=\linewidth]{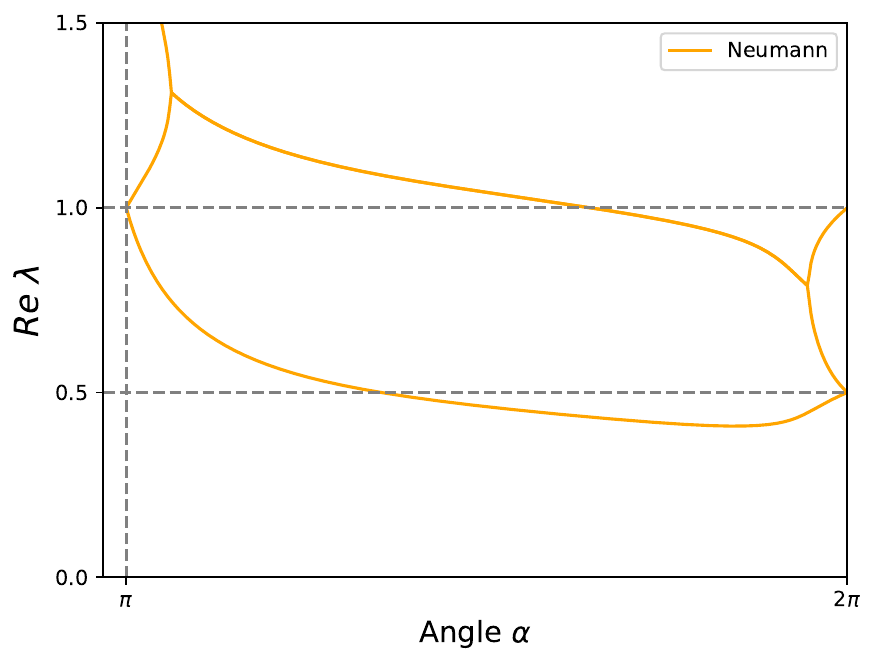}
    \end{subfigure}
    \hfill
    \begin{subfigure}[b]{0.48\linewidth}
        \centering
        \includegraphics[width=\linewidth]{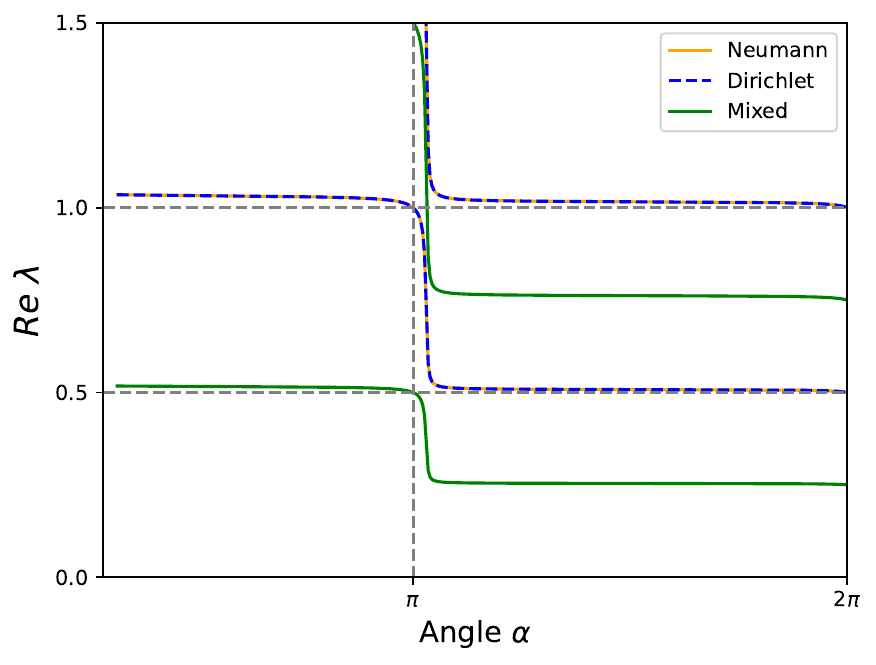}
    \end{subfigure}
    \caption{\small Left: Relation between $\Re \lambda$ and $\alpha\in [\pi,2\pi]$ for Neumann boundary conditions. The standard root of the monic elliptic tuple is given by $V=(S+i\Id_\ell)D$ for $S=\begin{pmatrix} 0&0\\0&2 \end{pmatrix}$, $D=\begin{pmatrix} 2&1\\1&2 \end{pmatrix}$. Although this tuple is contractive Neumann well-posed, it is not formal positive. The qualitative behavior differs from the Dirichlet and mixed case: Branch merging and $|\Re \lambda|<1/2$ occurs. Only selected branches are plotted. \\[1ex] Right: Relation between $\Re \lambda$ and $\alpha\in [1,2\pi]$ for different boundary conditions. The (scalar) elliptic tuple is defined by the standard root $V=-10+i$. In this case, the single branch for Neumann and Dirichlet boundary conditions coincides, and all branches approximate the bounds given in Theorem \ref{main dir} and \ref{main mixed}.}
    \label{fig:combined}
\end{figure}

We did not discuss bounds for $|\Re \lambda|$ in the case of Neumann boundary conditions. Recall that the corresponding condition in this case is $0\in\sigma(E~ Z_\alpha^{\lambda_+} E^{-1}- \overline{E} ~\overline{Z_\alpha}^{\lambda_-} \overline{E}^{-1})$, where $Z_\alpha$ is given in (\ref{Za}) and $E=\frac{1}{2}D^{-1/2}[S,D]D^{-1/2}+i\Id_\ell$. For $\alpha=\pi$, we have $Z^{\lambda_\pm}_\pi=\Id_\ell e^{\pm \lambda \pi}$, so the equation reduces to $\sin(\lambda \pi)=0$, as in the Dirichlet case (similarly for $\alpha=2\pi$). For all other angles, the bounds on $\Re \lambda$ are less clear. As Figure \ref{fig:combined} suggests, assuming contractive Neumann well-posedness does not guarantee $|\Re \lambda| >1/2$ for $\pi<\alpha<2\pi$. However, if the elliptic tuple $A$ is \textit{formal positive} (see Appendix \ref{ellli}), the literature (see §12 in \cite{VMR2001SPCS}) indicates that similar bounds to those for Dirichlet boundary conditions are obtained (for $\lambda \neq 0$).

\subsection{The scalar case}
In the scalar case ($\ell=1$), the matrices $A_\bullet$ reduce to real numbers, and the condition $0\in \sigma(M_{\lambda,\alpha})$ simplifies significantly. For Dirichlet and Neumann boundary conditions, the condition reduces to $Z_\alpha^{\lambda_+}=\overline{Z_\alpha}^{\lambda_-}$. For mixed boundary conditions, the condition reduces to $Z_\alpha^{\lambda_+}=-\overline{Z_\alpha}^{\lambda_-}$.

\subsection{Location of zeros}\label{localoca}
Theorem \ref{main dir} and Theorem \ref{main mixed} provide bounds on $|\Re \lambda|$, but they do not specify existence and location. In principle, Theorem \ref{main mixed} does not exclude the possibility of $\Lambda_\alpha=\emptyset$ for given $\alpha$. However, as Figure \ref{fig:1} suggests, this is not the case. For the particular case $\alpha \in \{\pi,2\pi\}$, it is straightforward to show that $1\in \Lambda_\pi$ and $1/2\in \Lambda_{2\pi}$ with multiplicity $\ell$ for Dirichlet and Neumann boundary conditions. This follows from $\sigma(\Id_\ell \sin(\lambda \pi))=\{\sin(\lambda \pi)\}$. Similarly, by using iii) in Theorem \ref{mix2}, one can show for mixed boundary conditions that there exist $\ell$ solutions $r^\lambda v$ (counted with multiplicity) satisfying $\Re \lambda =\frac{1}{2}$ for $\alpha=\pi$ ($\Re \lambda =\frac{1}{4}$ for $\alpha=2\pi$) if we assume contractive Neumann well-posedness. For angles $\alpha\notin\{\pi,2\pi\}$ existence of solutions in certain strips of the complex plane can be shown by using a generalization of Rouché's theorem. See the proof of Theorem 8.6.2 in \cite{VMR2001SPCS} for details.

\subsection{Optimality of the bounds}\label{herewediscuss}
The bounds on $|\Re \lambda|$ for different cases of $0<\alpha\leq 2\pi$ in Theorem \ref{main dir} and 1.) of Theorem \ref{main mixed} are sharp in the sense that one can construct a sequence of elliptic systems approaching these bounds. Note that the bounds for $\alpha\in \{\pi,2\pi\}$ are sharp, as discussed in Section \ref{localoca}. To illustrate this for other angles, it suffices to consider the scalar case. For $k\in \mathbb{N}$, define $S_k:=-k$ and $D_k:=1$. This yields the standard root $V_k=(S_k+i)D_k=-k+i$, and from (\ref{Za}) we derive:
\begin{align*}
    Z_{\alpha,k}=\cos(\alpha)-k\sin(\alpha)+i\sin(\alpha).
\end{align*}
For $0<\alpha<\pi$ and mixed boundary conditions, we obtain, as previously discussed, the condition $Z_{\alpha,k}^{\lambda}+\overline{Z_{\alpha,k}}^\lambda=0$ (here $\lambda_\pm=\lambda$). This is solved by $\lambda_{\alpha,k} \in \mathbb{R}$ such that $\Re(Z^{\lambda_{\alpha,k}}_{\alpha,k})=0$. As $\arg(Z_{\alpha,k})\xrightarrow[]{k\to \infty}\pi$, it follows that $\lambda_{k,\alpha}\xrightarrow[]{k\to \infty}\frac{1}{2}$ for any $0<\alpha<\pi$. Similar arguments can be given for $\pi<\alpha<2\pi$, and for Dirichlet and Neumann boundary conditions.

\subsection*{Contributions and Funding}
The author acknowledges the contributions of \emph{Joachim Rehberg}, who shared his expertise, provided further contacts in regularity theory, pointed out important references, and engaged in helpful discussions. The author acknowledges \emph{Matthias Liero} and \emph{Willem van Oosterhout} for feedback and fruitful discussions. The author acknowledges the hospitality of the University Kassel and helpful discussions with \textit{Dorothee Knees}. The author is funded by the Deutsche Forschungsgemeinschaft (DFG, German Research Foundation) under Germany's Excellence Strategy – The Berlin Mathematics
Research Center MATH+ (EXC-2046/1, project ID: 390685689).

\subsection*{Data Availability Statement}
 Some parts of this work, such as the plots, can be found as numerical implementation in the publicly available Jupyter Notebook: \href{https://doi.org/10.5281/zenodo.14417259}{https://doi.org/10.5281/zenodo.14417259}.

\appendix

\section{Ellipticity conditions}\label{ellli}
We aim to relate the ellipticity conditions presented in this work to those in the literature, in particular to ADN-elliptic systems \cite{ADN1964EBSE}.

\subsection*{Complementing boundary condition for $N_A\iff$ Neumann well-posedness}
The reference for the next paragraphs is §1.1.2 of \cite{MaRO2010EEPD}. Recall the setup given in \Cref{the model problem}: The domain $\mathcal{K}_{\alpha}$ with boundaries $\Gamma^\pm$ and the differential operator $L_{A}$ defined by $A=(A_{11},A_{12},A_{22})$ for $A_\bullet\in \Mat_\ell(\mathbb{R})$ symmetric and $A_{11},A_{22}>0$. For the boundary, we summarize the differential operators $B^{\pm}_{A}$ as $B_A(x)=B^{\pm}_{A}(x)$ for $x\in \Gamma^\pm$. Then, the system $(L_A,B_A)$ is called \textit{elliptic} (or ADN-elliptic) if the following two conditions are met:
\begin{enumerate}
    \item[1.)] The operator $L_A$ is \textit{properly elliptic}.
    \item[2.)] $B_A$ satisfies the \textit{complementing boundary condition} on $\partial \mathcal{K}_\alpha$.
\end{enumerate}
For 1.), in the case of real-valued $A_\bullet$, proper ellipticity of $L_A$ is equivalent to $L_A$ being elliptic (\ref{nostra}). For details, refer to §1 in \cite{Agmo1962EEGE}.

To address 2.), consider $x_0\in \partial \mathcal{K}_\alpha$ and $\xi$ tangential to $ \partial \mathcal{K}_\alpha$ at $x_0$. Denote by $\mathcal{M}(\xi)$ the subspace of solutions $u$ to the ODE:
\begin{align*}
    L_A(\xi-i n \partial_t) u(t)=0,\quad t>0,
\end{align*}
such that $u(t)\to 0$ for $t\to \infty$. Here, $n$ denotes the unit vector orthogonal to $\xi$. Now, $B_A$ is said to satisfy the complementing boundary condition if, for every $x_0\in \partial \mathcal{K}_\alpha$, every $\xi$ tangential to $\mathcal{K}_\alpha$ at $x_0$, and every $g\in \mathbb{C}^\ell$, there exists a unique $u\in \mathcal{M}(\xi)$ satisfying:
\begin{align*}
    B_A(\xi-in\partial_t ) u(t)|_{t=0}=g.
\end{align*}

Let us understand the complementing boundary condition for $\Gamma^-$ in the case of $B_A^-=N_A$. Then, $x_0\in \Gamma^-$, $\xi=(r,0)$ for $r\in \mathbb{R}\setminus \{0\}$, and $n=(0,1)$ such that
\begin{align*}
    L_A(\xi-i n\partial_t)&=r^2 A_{11}-2irA_{12}\partial_t-A_{22} \partial_t^2
\end{align*}
Let us assume $A$ is monic, and let $V$ be the standard root of $A$. Then for $r>0$ any $u\in \mathcal{M}(\xi)$ is of the form (recall Theorem \ref{algtheo} ii):
\begin{align*}
    u(t)=\exp(i rt V)c\quad \text{for }c\in \mathbb{C}^\ell\text{ and }t\in [0,\infty).
\end{align*}
Note that $\exp(irt \overline{V})$ does not occur, since $\|\exp(irt \overline{V})c\|\xrightarrow[]{t\to \infty} \infty$ for any $c\in \mathbb{C}^\ell\setminus \{0\}$ (deduced by the spectral mapping theorem). For $r<0$, we simply swap $V$ and $\overline{V}$. Continuing with $r>0$, we have
\begin{align*}
    N_A(\xi-in\partial_t)=\sum_{k,l=1}^2 A_{kl}n_k(\xi-in\partial_t)_l =r A_{12}-i\Id_\ell \partial_t\implies N_A(\xi-in\partial_t)u(t)\big|_{t=0}=r(A_{12}+V)c.
\end{align*}
So the complementing boundary condition reduces to the statement that, for any $g\in \mathbb{C}^\ell$, there exists a unique $c \in \mathbb{C}^\ell$ such that $(A_{12}+V)c =g$,
which is equivalent to invertibility of $A_{12}+V$. Comparing to Secion \ref{neumann boundary yea}, this shows that the complementing boundary condition for $B_A(x_0)=N_A(x_0)$ is equivalent to Neumann well-posedness. The arguments can also be generalized to nonmonic tuples and to $\Gamma^+$.

\subsection*{Formal positivity $\implies$ contractive Neumann well-posedness}
We have the following result for contractive Neumann well-posed systems.
\begin{lemma}\label{pathconnext}
Consider $A$ an elliptic tuple. Then:
\begin{enumerate}
    \item[i)] If $A$ is contractive Neumann well-posed, there exists a continuous path of contractive Neumann well-posed elliptic tuples connecting $A$ to $(\Id_\ell,0,\Id_\ell)$.
    \item[ii)] If there exists a continuous path from $A$ to $(\Id_\ell,0,\Id_\ell)$ consisting of Neumann well-posed tuples, then $A$ is contractive Neumann well-posed.
\end{enumerate}
\end{lemma}
\begin{proof}
i) Consider a contractive Neumann well-posed tuple $A=(A_{11},A_{12},A_{22})$. We can find a continuous path
\begin{align*}
        [0,1]\ni s\mapsto A(s)=(A_{11}(s),A_{12}(s),A_{22}(s)) \in (\Mat_\ell(\mathbb{R}))^3
\end{align*}
such that $A(1)=(\Id_\ell,0,\Id_\ell)$ and $A(0)=A$, with each $A(\bullet)$ being contractive Neumann well-posed. The path is constructed in three segments which can be glued together.\\
\textbf{First Segment: }Start with $(A_{11},A_{12},A_{22})$ and deform it as follows:
\begin{align*}
        [0,1]\ni s\mapsto A(s)=(A^{-s/2}_{22}A_{11}A^{-s/2}_{22},A^{-s/2}_{22}A_{12}A^{-s/2}_{22},A^{1-s}_{22}).
\end{align*}
Note that $A(0)=A$, $A(1)=(A^{-1/2}_{22}A_{11}A^{-1/2}_{22},A^{-1/2}_{22}A_{12}A^{-1/2}_{22},\Id_\ell)$, and that all $A(\bullet)$ are elliptic tuples and have the same monic reduction as $A(1)$. Thus, if $A(0)$ is contractive Neumann well-posed, then all $A(\bullet)$ are contractive Neumann well-posed.\\
\textbf{Second Segment: }Start with a contractive Neumann well-posed elliptic tuple of the form $ (A_{11},A_{12},\Id_\ell)$. Let us write $V=C+iD$ for its standard root and define $V_s:=(1-s)C+iD=((1-s)S+i\Id_\ell)D$ for $s\in [0,1]$. Using Theorem \ref{algtheo}, these represent the standard roots of the elliptic tuples:
\begin{align*}
        [0,1]\ni s\mapsto A(s)=\left(V_s^*V_s,-\frac{1}{2}(V_s^*+V_s),\Id_\ell\right).
\end{align*}
Note that $\rho([D^{-1},C])<2$ since $A(0)$ is contractive Neumann well-posed. So for $s\in [0,1]$:
    \begin{align*}
        \rho([(\Im V_s)^{-1},\Re V_s])=(1-s)\rho([D^{-1},C])< 2(1-s)<2.
\end{align*}
Moreover, we have $A(1)=(D^2,0,\Id_\ell)$ due to $-\frac{1}{2}(V^*_1+V^{}_1)=0$.\\
\textbf{Third Segment: }Start with an elliptic tuple of the form $(A_{11},0,\Id_\ell)$ and define $V_s:=iA_{11}^{1/2-s/2}$ for $s\in [0,1]$. By Theorem \ref{algtheo}, these are the standard roots of the elliptic tuples:
\begin{align*}
        [0,1]\ni s\mapsto A(s)=(A_{11}^{1-s},0,\Id_\ell).
\end{align*}
Note that $A(1)=(\Id_\ell,0,\Id_\ell)$ and $\rho([(\Im V_s)^{-1},\Re V_s])=0$.
\smallskip

ii) Assume there is a continuous path of elliptic tuples
\begin{align*}
        [0,1]\ni s\mapsto A(s)=(A_{11}(s),A_{12}(s),A_{22}(s)) \in (\Mat_\ell(\mathbb{R}))^3
\end{align*}
such that $A(1)=(\Id_\ell,0,\Id_\ell)$. Let $V_s=C_s+iD_s$ denote the standard root of the monic reduction of $A(s)$. Prove by contradiction and assume $A(0)$ is not contractive Neumann well-posed. Then there is $t>2$ such that $it\in \sigma([D_0^{-1},C_0])$ (recall Section \ref{neumann boundary yea}). Note that $\sigma([D_1^{-1},C_1])=\{0\}$, due to $C_1=0$ (see Example \ref{C=0,laplace}). So by continuity and $\sigma([D_s^{-1},C_s])\subset i \mathbb{R}$, there must be some intermediate value $s\in (0,1)$ such that $2i\in \sigma([D_s^{-1},C_s])$ and Neumann well-posedness is violated.
\end{proof}
A commonly used ellipticity condition in linear elasticity is \textit{formal positivity} (see §3.2 in \cite{CDN2010CSAR}), which is sometimes referred to as \textit{Legendre condition} (see §3.1.4 in \cite{GiMa2012IRTE}). The tuple $A=(A_{11},A_{12},A_{22})$ is said to be formal positive  if there exists $\kappa>0$ such that:
\begin{align}\label{shhimmer}
    \sum_{i,j=1}^2\langle A_{ij}f^{(i)},f^{(j)}\rangle \geq \kappa (\|f^{(1)}\|^2+\|f^{(2)}\|^2)\quad \text{for all } f^{(1)},f^{(2)}\in \mathbb{C}^\ell.
\end{align}
Note that the LHS of (\ref{shhimmer}) can be rewritten in block matrix form as $f^* M_A f\geq \kappa \|f\|^2$, where
\begin{align*}
    f:=\begin{pmatrix}
        f^{(1)} & f^{(2)}
    \end{pmatrix}\in \mathbb{C}^{2\ell},\quad 
    M_A:=\begin{pmatrix}
        A_{11}&A_{12}\\A_{12}&A_{22}
    \end{pmatrix}\in \Mat_{2\ell}(\mathbb{R}).
\end{align*}
Thus, formal positivity is equivalent to $M_A>0$. Next, we state Theorem 7.7.6 in \cite{Horn1985MaAn}:
\begin{theorem}\label{blocki}
    Consider $A,B,C\in \Mat_\ell(\mathbb{C})$ and the block matrix $
        M=\begin{pmatrix}
        A&B\\B^*&C
    \end{pmatrix}$. We have $M>0$ if and only if $A>0$ and $C>B^*A^{-1}B$. Furthermore, this is equivalent to $\rho(B^*A^{-1}BC^{-1})<1$.
\end{theorem}

The next result relates formal positivity and contractive Neumann well-posedness.
\begin{lemma}\label{seemeti}
    Consider an elliptic tuple $A$. If it is formal positive, then it is contractive Neumann well-posed.
\end{lemma}
\begin{proof}
Consider an elliptic tuple $A$. Note that $A$ is formal positive if and only if its monic reduction is formal positive. To see this, apply $\operatorname{diag}(A_{22}^{-1/2},A_{22}^{-1/2})$ from the left and right to $M_A$ and use Sylvester's law of inertia. So, without loss of generality, assume $A_{22}=\Id_\ell$.\smallskip

\textbf{Claim: }If $A$ is formal positive, then $A$ is Neumann well-posed.
\smallskip 

Before proving the claim, let us argue why the claim shows the statement. Due to $M_A>0$ and Theorem \ref{blocki}, we have:
\begin{align*}
        M_{A}(s):=\begin{pmatrix}
        A_{11}&sA_{12}\\s A_{12}&\Id_\ell
    \end{pmatrix}>0\quad \text{for all } s\in [0,1].
\end{align*}
Thus, using the claim, any $A(\bullet)$ on the path $s\mapsto A(s)=(A_{11},s A_{12},\Id_\ell)$ is Neumann well-posed. Using techniques from the proof of Lemma \ref{pathconnext}, we can further deform $(A_{11},0,\Id_\ell)$ to $(\Id_\ell,0,\Id_\ell)$, showing that $A$ is path-connected to $(\Id_\ell,0,\Id_\ell)$ by a path of Neumann well-posed systems. From Lemma \ref{pathconnext}, it follows that $A$ is contractive Neumann well-posed. It remains to prove the claim.
\smallskip

\textit{Proof of the claim: }Using Theorem \ref{algtheo}, we write for the standard root $V=(S+i\Id_\ell)D$ for $S,D\in \Mat_\ell(\mathbb{R})$ symmetric and $D>0$, and $M_A$ as:
\begin{align*}
        M_A=\begin{pmatrix}
            D(S^2+\Id_\ell)D&-\frac{1}{2}(SD+DS)\\
            -\frac{1}{2}(SD+DS)&\Id_\ell
        \end{pmatrix}.
\end{align*}
Now, assume $M_A>0$ and $A$ not Neumann well-posed, and derive a contradiction. The latter implies $2i\in \sigma(D^{-1}SD-S)$ (recall (\ref{recalli})). $M_A>0$ is by Theorem \ref{blocki} equivalent to:
\begin{align}\label{cannoot}
       \rho((SD+DS)(D(S^2+\Id_\ell)D)^{-1}(SD+DS))<4\iff \rho(N^*N)<4\iff \|N\|<2,
\end{align}
where $N:=(S+i\Id_\ell)^{-1}(D^{-1}SD+S)$. By assumption, there exists $y\in \mathbb{C}^\ell$ with $\|y\|=1$ such that $(D^{-1}SD-S)y=2i y$. This implies $(D^{-1}SD+S)y=2(S+i\Id_\ell) y$, and thus $ \|Ny\|=2$, contradicting (\ref{cannoot}). This proves the claim.
\end{proof}

\section{Functional calculus}\label{spundfunc}
We briefly summarize the functional calculus adapted to finite-dimensional vector spaces. The reference is \textit{Symbolic Calculus} in §10 of \cite{Rudin1991FA}. Consider a holomorphic function $f:\Omega \subset \mathbb{C}\to \mathbb{C}$ on an open set $\Omega $ and let $A\in \Mat_\ell(\mathbb{C})$ with $\sigma(A)\subset \Omega $. We define
    \begin{align}\label{funcca}
        f(A):=\frac{1}{2\pi i} \int_\Gamma f(z) (\Id_\ell \cdot z-A)^{-1}dz,
    \end{align}
where $\Gamma$ is a contour enclosing $\sigma(A)$. The definition is independent of the choice of contour $\Gamma$, provided it encloses $\sigma(A)$. Then, $(\Id_\ell \cdot z-A)^{-1}$ and the integral are well-defined. The functional calculus exhibits the following properties:
\begin{itemize}
    \item \textit{Spectral mapping theorem}: The spectrum of $f(A)$ satisfies $ \sigma(f(A))=f(\sigma(A))$.
    \item For any invertible $Q\in \Mat_\ell(\mathbb{C})$, it holds that $f(QAQ^{-1})=Qf(A)Q^{-1}$.
\end{itemize}

\subsection*{Complex Exponentiation}
For $\lambda \in \mathbb{C}$ and $z\in \mathbb{C}\setminus \{0\}$, complex exponentiation was defined in Section \ref{explisolf} by $z^{\lambda_a}=\exp(\lambda \log_a(z))$. Using the functional calculus, this is generalized to $A^{\lambda_a}$ for $A\in \Mat_\ell(\mathbb{C})$ provided $0\notin \sigma(A)$. The following exponential rules hold for $\lambda,\mu\in \mathbb{C}$ and $a\in \{+,-,o\}$:
\begin{align*}
    A^{(\lambda+\mu)_a}=A^{\lambda_a}A^{\mu_a},\quad (A^{\lambda_a})^{\mu_a}=A^{(\lambda\cdot \mu)_a},\quad \overline{A^{\lambda_a}}=\overline{A}^{\overline{\lambda}_a}.
\end{align*}
Also, by the spectral mapping theorem, $0\notin \sigma(A^{\lambda_a})$. Most of the time, we consider $Z^{\lambda_a}$ for $Z\in \Mat_\ell(\mathbb{C})_{+i}$ (see Def. \ref{def+i}). This is well-defined due to $\sigma(Z)\cap \mathbb{R}=\emptyset$ (see Lemma \ref{symm}). Also, $(Z^{\lambda_a})^T=Z^{\lambda_a}$. This follows from the Cauchy integral formula and the fact that $\Id_\ell z-Z$, and thus its inverse, are symmetric matrices for any $z\in \mathbb{C}$.

\section{Accretive operators}\label{numerrange and accreive}
The reference for this section is the dissertation \cite{Haas2003FCSO}. Throughout this part, the $\arg$- and $\log$-function correspond to the principal branch. The goal is to prove Lemma \ref{numrang} and Lemma \ref{chaotisch}.\smallskip

The class of matrices $\Mat_\ell(\mathbb{C})_{+i}$ is naturally closely related to matrices $\Mat_\ell(\mathbb{C})\ni A=A^T$ with $\Re A>0$, which are examples of \textit{accretive operators}. For $0<\omega\leq \pi$, we define:
\begin{align*}
    S_\omega:=\{z\in \mathbb{C}: z\neq 0\text{ and }|\arg(z)|<\omega \}.
\end{align*}

\begin{definition}[p.101 in \cite{Haas2003FCSO}]
    Let $0\leq \omega\leq \frac{\pi}{2}$. $A\in \Mat_\ell(\mathbb{C})$ is called $\omega$-\textit{accretive} if $W(A)\subset \operatorname{clos}(S_\omega)$. For $\omega=\frac{\pi}{2}$, i.e., $W(A)\subset \operatorname{clos}(\operatorname{RHS})$, $A$ is simply called \textit{accretive}.
\end{definition}

\begin{theorem}[Theorem B.21 in \cite{Haas2003FCSO}]\label{kontraka}
    $A\in \Mat_\ell(\mathbb{C})$ is accretive if and only if $-A$ generates a strongly continuous contraction semigroup, i.e., $ \|e^{-At}\|\leq 1$ for all $t\geq 0$.
\end{theorem}
The numerical range of fractional powers of accretive operators is well understood. The following results, all found in §\textit{Fractional Powers of m-Accretive Operators and the Square Root Problem} of \cite{Haas2003FCSO}, will be used to extend these results to $\Mat_\ell(\mathbb{C})_{+i}$.

\begin{proposition}\label{amack}
    Let $\delta>0$ and $A-\delta\Id_\ell $ be accretive for $A\in \Mat_\ell(\mathbb{C})$. Then $A^\lambda-\delta^\lambda \Id_\ell$ is accretive for each $0<\lambda\leq 1$.
\end{proposition}

\begin{proposition}\label{wingel}
    Let $A\in \Mat_\ell(\mathbb{C})$ be accretive, and let $0\leq  \lambda\leq  1$. Then $A^\lambda $ is $\frac{\lambda \pi}{2}$-accretive, i.e., $W(A^\lambda)\subset \operatorname{clos}(S_{\frac{\lambda \pi}{2}})$.
\end{proposition}

\begin{proposition}\label{loggi}
    Let $A\in \Mat_\ell(\mathbb{C})$ be an injective $\omega$-accretive operator for some $0\leq \omega \leq \frac{\pi}{2}$. Then $W(\log(A))\subset \{z\in \mathbb{C}: |\Im z|\leq \omega \}$.
\end{proposition}

We need two preparatory Lemmas.
\begin{lemma}\label{lachflash}
    Consider $0<\lambda \leq 1$ and $A\in \Mat_\ell(\mathbb{C})$ with $W(A)\subset \operatorname{RHP}$. Then $W(A^\lambda)\subset S_{\frac{\lambda\pi}{2}}$.
\end{lemma}
\begin{proof}\label{bursti}
    Since $W(A)\subset \operatorname{RHP}$, by \textbf{N2} there exists some $\delta>0$ such that $A-\delta \Id_\ell$ is accretive. Applying Prop. \ref{amack}, we deduce $W(A^\lambda-\delta^\lambda\Id_\ell)\subset \operatorname{clos}(S_{\frac{\lambda\pi}{2}})$ and the claim follows by \textbf{N2}.
\end{proof}
\begin{lemma}\label{uspis}
    Consider $Z\in \Mat_\ell(\mathbb{C})_{+i}$. Then $-\Im (Z^{-1})>0$.
\end{lemma}
\begin{proof}
    Note that 
    \begin{align*}
        -\Im (Z^{-1})=\frac{1}{2i}(\overline{Z}^{-1}-Z^{-1})=\frac{1}{2i}\overline{Z}^{-1}(Z-\overline{Z})Z^{-1}=(Z^{-1})^* \Im (Z) Z^{-1}>0,
    \end{align*}
    where the last line follows from $\Im Z>0$ and Sylvester's law of inertia.
\end{proof}

Now that we established the necessary tools, we prove Lemma \ref{numrang} and Lemma \ref{chaotisch}.

\begin{proof}[Proof of Lemma \ref{numrang}]
Consider $Z\in \Mat_\ell(\mathbb{C})_{+i}$ and $0<\lambda\leq 1$, and show $W'(Z^\lambda)\subset \operatorname{UHP}$. Define $A=- i Z$ so $\Re A=\Im Z>0$, implying $W(A )\subset \operatorname{RHP}$. Using this, we bound the numerical range of $Z^\lambda$ as follows:
\begin{align}\label{noteitnote}
    W(Z^\lambda)=&~i^\lambda W(A^\lambda)=e^{i\lambda \pi/2}W(A^\lambda)\subset e^{ i\lambda \pi /2} S_{\frac{\lambda\pi}{2}}\subset \operatorname{UHP},
\end{align}
where we applied \textbf{N2} and Lemma \ref{lachflash}. Note that this implies $W'(Z^\lambda)\subset \operatorname{UHP}$ as well. Next, assume $-1\leq \lambda<0$, and show $W'(Z^\lambda)\subset -\operatorname{UHP}$. Note that $\overline{Z^{-1}}\in \Mat_\ell(\mathbb{C})_{+i}$ by Lemma \ref{uspis}. The result is deduced by
\begin{align*}
    \overline{W(Z^\lambda)}\stackrel{(\textbf{N8})}{=}W(\overline{Z^\lambda})=W(\overline{Z^{-1}}^{|\lambda|})\stackrel{(\ref{noteitnote})}{\subset}\operatorname{UHP}\implies W(Z^\lambda)\subset \overline{\operatorname{UHP}}=-\operatorname{UHP}.
\end{align*}
\end{proof}

\begin{proof}[Proof of Lemma \ref{chaotisch}] Consider $Z\in \Mat_\ell(\mathbb{C})_{+i}$ and $\lambda \in \mathbb{R}\setminus \{0\}$.
\smallskip

\textbf{Claim 1: }$\|Z^{i\lambda}\|<1$ for $\lambda>0$.

Note that $Z^{i\lambda}=\exp(i\lambda \log(Z))$, and due to Theorem \ref{kontraka}, it suffices to show that $-iW(\log(Z))\subset \operatorname{RHP}$. Let $A=-i Z$, so $W(A)\subset \operatorname{RHP}$ and $A$ is $\frac{\pi-\varepsilon}{2}$-accretive for some $\varepsilon>0$. Claim 1 is deduced by
\begin{align*}
    -iW(\log(Z))=&~-iW(\log(iA))=-iW(\log(i)\Id_\ell+\log(A))=-i W\Big(\frac{i\pi}{2}\Id_\ell\Big)-i W(\log(A))\\
    \subset &~\frac{\pi}{2}-i\cdot \left\{z\in \mathbb{C}: |\Im z|\leq \frac{\pi-\varepsilon}{2} \right\}\subset \operatorname{RHP},
\end{align*}
where we used properties of the $\log$-function, functional calculus, \textbf{N5}, \textbf{N2}, and Prop. \ref{loggi}.
\smallskip

\textbf{Claim 2: }$\rho((Z^{i\lambda })^*Z^{i\lambda})\xrightarrow[]{\lambda\to \infty} 0$.

For any $N\in \mathbb{N}$, we have:
\begin{align*}
    \rho((Z^{i\lambda })^*Z^{i\lambda})=\|(Z^{i\lambda })^*Z^{i\lambda})\|=\|Z^{i\lambda}\|^2=\|(Z^{i\lambda/N})^N\|^2\leq \|Z^{i\lambda/N}\|^{2N},
\end{align*}
where we used $\|B^N\|\leq \|B\|^N$. Claim 1 implies $ \rho((Z^{i\lambda })^*Z^{i\lambda})\xrightarrow[]{\lambda\to \infty} 0$.
\smallskip

i) Proof of $\operatorname{sgn}(\lambda)\left(\Id_\ell-(Z^{i\lambda })^*Z^{i\lambda}\right)>0$: Note that for any $A\in \Mat_\ell(\mathbb{C})$ one has $\|AA^*\|=\|A^*A\|=\|A\|^2=\rho(A^*A)$, so $\Id_\ell-A^*A>0$ is equivalent to $\|A\|<1$. Thus,  the case $\operatorname{sgn}(\lambda)>0$ follows from Claim 1. The case $\operatorname{sgn}(\lambda)<0$ can be derived from the case $\operatorname{sgn}(\lambda)>0$, the fact 
\begin{align*}
    \Id_\ell-B>0 \iff B^{-1}-\Id_\ell>0\quad \forall B\in \Mat_\ell(\mathbb{C}) \text{ with }B=B^*
\end{align*}
applied to $B=(Z^{i\lambda })^*Z^{i\lambda}$, and $\sigma((Z^{-i\lambda })^*Z^{-i\lambda})=\sigma(Z^{-i\lambda} (Z^{-i\lambda })^*)$.

ii) Proof of $\min\{\beta: \beta \in \sigma((Z^{i\lambda })^*Z^{i\lambda})\}\xrightarrow[]{\lambda\to -\infty} \infty$: This follows from Claim 2 using $(Z^{i\lambda})^{-1}=Z^{-i\lambda}$, $\sigma((Z^{-i\lambda })^*Z^{-i\lambda})=\sigma(Z^{-i\lambda} (Z^{-i\lambda })^*)$, and the spectral mapping theorem.
\end{proof}

\end{document}